\documentclass[12pt,twosided]{article}

\usepackage{amsmath}
\usepackage{amssymb}
\usepackage{theorem}
\usepackage{epsfig}
\usepackage{xypic}
\usepackage{latexsym}
\usepackage{rotating}
\usepackage{times}

\def\sym{\mathbb}
\def\gl{{\mathfrak gl}}

\def\F{\mathbf{F}}

\def\G{\mathbf{G}}

\def\C{{\sym C}}

\def\Z{{\sym Z}}
\def\Q{{\sym Q}}
\def\L{{\cal L}}
\def\SG{{\mathfrak S}}

\def\std{\operatorname{Std}}
\def\Des{\operatorname{Des}}
\def\Rec{\operatorname{Rec}}
\def\SS{\mathbf{S}}
\def\Sh{\operatorname{Sh}}
\def\FQSym{\mathbf{FQSym}}
\def\MQSym{\mathbf{MQSym}}
\def\FSym{\mathbf{FSym}}
\def\Sym{\mathbf{Sym}}
\def\MS{\mathbf{MS}}
\def\FS{{\sym \phi S}}
\def\T{\mathbf{T}}

\def\P{\mathbf{P}}
\def\HQ{{\sf Q}}
\def\HR{{\sf R}}
\def\dim{\operatorname{dim}}
\def\rad{\operatorname{rad}}
\def\soc{\operatorname{soc}}
\def\mod{\operatorname{mod}}
\def\Ext{\operatorname{Ext}}
\def\conn{{\cal C}}
\def\ff{{\sf F}}
\def\V{\mathbf{V}}
\def\sfact{\operatorname{\rm sfact}}
{
\theoremstyle{plain}
\newtheorem{theorem}{Theorem}[section]

\newtheorem{corollary}[theorem]{Corollary}
\newtheorem{conjecture}[theorem]{Conjecture}
\newtheorem{definition}[theorem]{Definition}
\newtheorem{proposition}[theorem]{Proposition}

\newtheorem{lemma}[theorem]{Lemma}

\theorembodyfont{\rmfamily}
\newtheorem{example}[theorem]{Example}

\newtheorem{remark}[theorem]{Remark}
}

\def\Proof{\noindent \it Proof -- \rm}
\def\cqfd{\hspace{3.5mm} \hfill \vbox{\hrule height 3pt depth 2 pt width 2mm}
\bigskip}
\def\qed{\hspace{3.5mm} \hfill \vbox{\hrule height 3pt depth 2 pt width 2mm}}

\def\U0{{\cal U}_0(gl_N)}

\newcommand{\free}[1]{\langle#1\rangle}
\newcommand{\pair}[2]{\left\langle\,#1\ \vrule\ #2\,\right\rangle}
\newcommand{\pairing}[2]{\left\langle#1,#2\right\rangle} 

\def\prob{\operatorname{Prob}}
\def\maj{\operatorname{maj}}

\def\ie{{\it i.e.\ }}
\def\cf{{\it cf.\ }}

\def\MM{\mathbf{M}}
\def\NN{\mathbf{N}}
\def\ch{\operatorname{ch}}
\def\fch{\mathbf{ch}}
\def\K{{\sym K}}

\def\bGamma{{\bf \Gamma}}

\def\shuff#1#2{\mathbin{
      \hbox{\vbox{
        \hbox{\vrule
              \hskip#2
              \vrule height#1 width 0pt
               }%
        \hrule}%
             \vbox{
        \hbox{\vrule
              \hskip#2
              \vrule height#1 width 0pt
               \vrule }%
        \hrule}%
}}}
\def\shuffl{{\mathchoice{\shuff{7pt}{3.5pt}}%
                        {\shuff{6pt}{3pt}}%
                        {\shuff{4pt}{2pt}}%
                        {\shuff{3pt}{1.5pt}}}}%
\def\shuffle{\, \shuffl \,}

\def\ash{\, \underline{\shuffl} \, }

\def\DessinMatrix#1{\vcenter{\hbox{\makebox[1.7ex]{$#1$}}}}
\def\GenMatrix#1{\vcenter{\halign{&$\DessinMatrix{##}$\cr#1}}\egroup}

\def\setinterlineskip#1{\baselineskip=0pt
  \lineskip=#1 \lineskiplimit=\maxdimen}

\def\matrice{%
  \bgroup
  \let\ =\omit
  \let\\=\cr
  \setinterlineskip{4.0pt}
  \GenMatrix}

\def\DessinsMatrix#1{\vcenter{\hbox{\makebox[1.3ex]{$\scriptstyle#1$}}}}
\def\GensMatrix#1{\vcenter{\halign{&$\DessinsMatrix{##}$\cr#1}}\egroup}

\def\smallmatrice{%
  \bgroup
  \let\ =\omit
  \let\\=\cr
  \setinterlineskip{3.0pt}
  \GensMatrix}

\def\tensor{\otimes}

\newcommand{\SMat}[1]{{\left[\smallmatrice{#1\\}\right]}}

\newlength{\Hackl}
\newcommand{\Hack}{\vrule height \Hackl width 0pt}
\newcommand{\indexmat}%
    {\smallmatrice{\Hack a\\\Hack b\\\Hack c\\\Hack d\\\Hack e\\}}

\newdimen\Squaresize \Squaresize=14pt
\newdimen\Thickness \Thickness=0.5pt

\def\Square#1{\hbox{\vrule width \Thickness
   \vbox to \Squaresize{\hrule height \Thickness\vss
      \hbox to \Squaresize{\hss#1\hss}
   \vss\hrule height\Thickness}
\unskip\vrule width \Thickness}
\kern-\Thickness}

\def\Vsquare#1{\vbox{\Square{$#1$}}\kern-\Thickness}

\def\young#1{
\vbox{\smallskip\offinterlineskip
\halign{&\Vsquare{##}\cr #1}}}

\title{\LARGE \sf 
Noncommutative Symmetric Functions VI: \\
Free Quasi-Symmetric Functions \\
and Related Algebras}

\author{G\'erard~{\sc Duchamp},
        Florent~{\sc Hivert} and
        Jean-Yves~{\sc Thibon}
}

\date{}

\begin{document}
\thispagestyle{empty}
\maketitle

\begin{abstract}
This article is devoted to the study of several algebras which are related  to symmetric functions,
and which admit  linear bases labelled by various combinatorial objects:
permutations (free quasi-symmetric functions), standard Young tableaux
(free symmetric functions) and packed integer matrices (matrix quasi-symmetric
functions). 
Free quasi-symmetric functions provide a kind of noncommutative
Frobenius characteristic for a certain category of modules over the $0$-Hecke algebras. 
New examples of indecomposable $H_n(0)$-modules are discussed, and the homological
properties of $H_n(0)$ are computed for small $n$. Finally, the algebra of matrix quasi-symmetric functions
is interpreted as a convolution algebra.
\end{abstract}

\newpage
\tableofcontents

\newpage

                       \section{Introduction}

This article is devoted to the study of several algebras closely related
to symmetric functions. By `closely related', we mean that these
algebras can be fitted into a diagram of homomorphisms
$$
\newdir{ >}{{}*!/-4mm/@{>}}
\xymatrix@R=10mm@L=10mm{
{\FSym}\ar@{ >->}[rr]   
\ar@{->>}[dd]     &               &{\FQSym}\ar@{ >->}[dr]
                                                   \ar@{->>}[d]     &  \\
                & {\Sym}\ar@{ >->}[ul]
                        \ar@{->>}[dl]
                        \ar@{ >->}[ur] 
                        \ar[dr]   &QSym_q\ar@{->>}[d] &{\MQSym}
                                                           \ar@{->>}[dl]\\
{Sym}\ar@{ >->}[rr] &                    &Qsym               &  \\
}
$$
along which most of the interesting structure can pulled back or
pushed forward.

Our notation is summarized in the following table, which indicates
also the combinatorial objects labelling the natural bases of the
various algebras:

\bigskip

\begin{center}
{\footnotesize
\begin{tabular}{|c|c|c|}
\hline
Symbol & Algebra & Basis \\
\hline\hline
$Sym$ & Symmetric functions & Partitions \\
\hline
$Qsym$ & Quasi-symmetric functions & Compositions\\
\hline
$\Sym$ & Noncommutative symmetric functions & Compositions \\
\hline
$QSym_q$ & Quantum quasi-symmetric functions & Compositions ($\C(q)$-basis)\\
\hline
$\FSym$ & Free symmetric functions & Standard Young tableaux \\
\hline
$\FQSym$ & Free quasi-symmetric functions & Permutations \\
\hline
$\MQSym$ & Matrix quasi-symmetric functions & Packed integer matrices \\
\hline
\end{tabular}
}
\end{center}
\bigskip

The starting point of the construction is the triangular diagram
formed by the embedding of $Sym$ in $QSym$, and the abelianization map
from $\Sym$ onto $Sym$ \cite{NCSF1}.  Since the maps preserve the
natural gradations, we have, for the homogeneous components of degree
$n$ of these algebras, a commutative diagram
$$
\xymatrix{
Sym_n \ar[rr]^d &                             & QSym_n \\
                & {\Sym_n}\ar[lu]_e \ar[ru]_c & \\
}
$$
which has a neat interpretation in representation theory: it is the
Cartan-Brauer triangle of $H_n(0)$, the Hecke algebra of type
$A_{n-1}$ at $v=0$ \cite{DKLT,NCSF3,NCSF4}. This means that $QSym_n$
is to be be interpreted as $G_0(H_n(0))$, the Grothendieck group of
the category of finitely generated $H_n(0)$-modules, $\Sym_n$ as
$K_0(H_n(0))$, the Grothendieck group of finitely generated projective
$H_n(0)$-modules, and $Sym_n$ as the Grothendieck group
$R(H_n(v))=G_0(H_n(v))=K_0(H_n(v))$ of the semi-simple algebra
$H_n(v)$, for generic $v$. Moreover, the inclusion map $d:\ 
Sym_n\rightarrow QSym_n$ is the decomposition map. Indeed, the simple
$H_n(v)$ modules $V_\lambda(v)$ correspond to the Schur functions
$s_\lambda$, with $\lambda\vdash n$, and the coefficients of the
quasi-symmetric expansion
\begin{equation}
s_\lambda=\sum_{|I|=n}d_{\lambda I}F_I
\end{equation}
are the multiplicities of the simple $H_n(0)$ modules $\SS_I$
(parametrized by compositions $I$ of $n$) as composition factors of
the specialized module $V_\lambda(0)$.

This interpretation leads to a $q$-analogue of $QSym$: the algebra
$QSym_q$ of quantum quasi-symmetric functions, defined in \cite{TU}.
Here, the indeterminate $q$ is introduced to record a certain
filtration on $H_n(0)$-modules. For generic complex values of $q$,
$QSym_q$ is non commutative, and in fact isomorphic to $\Sym$, but for
$q=1$ one recovers the commutative algebra of quasi-symmetric functions
$QSym$.
\smallskip

This construction can be somewhat clarified by the introduction of the
larger algebra $\FQSym$, a subalgebra of the free associative algebra
$\C\free{A}$ (whence the name free quasi-symmetric functions) which
admits $\Sym$ as a subalgebra, and is mapped onto $QSym_q$ when one
imposes the $q$-commutation relations of the quantum affine space
($a_ja_i=qa_ia_j$ for $j>i$) on the letters of $A$.

This algebra turns out to be isomorphic to the convolution algebra of
symmetric groups studied by Malvenuto and Reutenauer \cite{MR1}. It
contains a subalgebra whose bases are naturally labelled by standard
Young tableaux, which provides a concrete realization of the algebras
of tableaux of Poirier and Reutenauer \cite{PR}. We call it $\FSym$,
the algebra of free symmetric functions. To illustrate the relevance
of the realization of $\FSym$ as an algebra of noncommutative
polynomials, we use it to present a complete proof of the
Littlewood-Richardson rule within a dozen of lines (the idea of the
proof is not new, but the formalism makes it quite compact and
transparent).  In the same vein, we show that the use of $\FQSym$
allows one to give simple presentations of Stanley's QS-distribution
\cite{St} and of the Hopf algebra of planar binary trees of Loday and
Ronco \cite{LR}.
\medskip

The next step is to look for a representation theoretical
interpretation of $\FQSym$. It turns out that $\FQSym_n$ can be
interpreted as a kind of Grothendieck group for a certain category
${\cal N}_n$ of $H_n(0)$-modules,  which contains in particular
simple, projective, and skew Specht modules. However, this is far from
exhausting all the $H_n(0)$-modules, since we prove that for $n\ge 4$,
$H_n(0)$ is not representation finite. As a step towards a more
exhaustive study of the $0$-Hecke algebras, we determine their quivers
for all $n$, and discuss their homological properties for small values
of  $n$.
\medskip

Finally, we show that $\FQSym$ can be embedded into a larger algebra,
$\MQSym$, whose bases are labelled by packed integer matrices, or, if
one prefers, by double cosets of symmetric groups modulo parabolic
subgroups. This is a self-dual bialgebra, which accommodates all the
previous ones as quotients or subalgebras, and in which most of the
structure of symmetric functions survives. It is not known whether
$\MQSym$ can be interpreted as a sum of Grothendieck groups. It has,
however, some representation theoretical meaning, as the centralizer
algebra of $GL(N,\C)$ in a certain infinite dimensional
representation.

\medskip
{\it Acknowledgements.} This paper was completed during the stay
of the authors at the Isaac Newton Institute for Mathematical Science,
whose hospitality is gratefully acknowledged. F. H. and G. D. were
supported by the European Community, and J.-Y. T. by an EPSRC grant.

\pagebreak[4]
                      \section{Background}

\subsection{Hypoplactic combinatorics}

Our notations will be essentially as in \cite{NCSF1}. In this paper,
we will use the realization of $\Sym$ as a subalgebra of the free
associative algebra $\C\free{A}$ over an infinite ordered noncommutative
alphabet $A=\{a_i\ | \ i\ge 1\}$. Then, the ribbon Schur function
$R_I$ is identified with the sum of all words whose shape is encoded
by the composition $I$.

We recall the notion of quasi-ribbon words and tableaux. A
quasi-ribbon tableau of shape $I$ is a ribbon diagram $r$ of sha\-pe
$I$ filled by letters of $A$ in such a way that each row of $r$ is
nondecreasing from left to right, and each column of $r$ is strictly
increasing from top to bottom. A word is said to be a quasi-ribbon
word of shape $I$ if it can be obtained by reading from bottom to top
and from left to right the columns of a quasi-ribbon diagram of shape
$I$.
	                                
The hypoplactic Robinson-Schensted correspondence is a bijection
between words $w$ and pairs $(\HQ(w),\HR(w))$, where $\HQ(w)$ and
$\HR(w)$ are respectively a quasi-ribbon tableau and a standard ribbon
tableau of the same shape \cite{NCSF4}. The equivalence relation on
words $u$ and $v$ defined by
\begin{equation}
  u\equiv v \ \Longleftrightarrow \HQ(u)  =  \HQ(v)
\end{equation}
can be shown to coincide with the hypoplactic congruence of the free
monoid $A^*$, which is generated by the plactic relations
\begin{equation*}
  \left\{ \,
    \begin{array}{c@{\ }ccl}
      aba \, \equiv  \, baa & , &
      bba \, \equiv  \, bab & \qquad \hbox{for \ $a<b$,} \\[2mm]
      acb \, \equiv  \, cab & , &
      bca \, \equiv  \, bac & \qquad \hbox{for \ $a<b<c$}.
    \end{array}
  \right.
\end{equation*}
and the quartic hypoplactic relations
\begin{equation*}
  \left\{ \,
    \begin{array}{c@{\ }ccl}
      baba \, \equiv  \, abab & , & baca \, \equiv  \, abac & \qquad
                                       \text{for \ $a<b<c$,} \\[2mm]
      cacb \, \equiv  \, acbc & , &
      cbab \, \equiv  \, bacb & \qquad \text{for \ $a<b<c$,} \\[2mm]
      badc \, \equiv  \, dbca & , &
      acbd \, \equiv  \, cdab & \qquad \text{for \ $a<b<c<d$}.
    \end{array}
  \right.
\end{equation*}
Despite the apparent complexity of these relations, it can be shown
that $u\equiv v$ if and only if $u$ and $v$ have the same evaluation
and the permutations $\std(u)^{-1}$ and $\std(v)^{-1}$ have the same
descents. Here, $\std(w)$ denotes the standardized of the word $w$,
{\it i.e.} the permutation obtained by iteratively scanning $w$ from
left to right, and labelling $1,2,\ldots$ the occurrences of its
smallest letter, then numbering the occurrences of the next one, and
so on. Alternatively, $\sigma=\std(w)^{-1}$ can be characterized as
the unique permutation of minimal length such that $w\sigma$ is a
nondecreasing word.

Quasi-symmetric functions can be lifted to the hypoplactic algebra.
The hypoplactic quasi-ribbon $F_I(A)$ is defined as the sum of all
quasi-ribbon words of shape $I$ in the hypoplactic algebra. It is
shown in \cite{NCSF4} that these elements span a commutative
$\Z$-subalgebra, and that the image of $F_I(A)$ in $\Z[X]$ by the
natural homomorphism is the usual quasi-symmetric function $F_I(X)$.

\subsection{$0$-Hecke algebras}

The $0$-Hecke algebra $H_n(0)$ is the $\C$-algebra generated by $n-1$
elements $T_1,\ldots, T_{n-1}$ satisfying the braid relations and
$T_i^2=-T_i$. It will be convenient to introduce a special notation
for the generators $\xi_i=1+T_i$ and $\eta_i=-T_i$, which also satisfy
the braid relations, $\xi_i^2=\xi_i$ and $\eta_i^2=\eta_i$. To a
permutation $\sigma\in\SG_n$, we can therefore associate three
elements $T_\sigma$, $\xi_\sigma$ and $\eta_\sigma$ by the usual
process of taking the products of generators labelled by a reduced
word for $\sigma$.

The irreducible $H_n(0)$ modules are denoted by $\SS_I$ and the unique
indecomposable projective module $M$ such that $M/\rad (M)=\SS_I$ is
denoted by $\P_I$. Its socle is simple and isomorphic to
$\SS_{\bar I}$, where $\bar I$ is the mirror composition of $I$.

The dimension of $\P_I$ is equal to the cardinality of the descent
class $D_I$, the set of permutations having $I$ as descent
composition. This set is an interval $[\alpha(I),\omega(I)]$ of the
(left) weak order on $\SG_n$. As shown by Norton \cite{No}, one can
realize $\P_I$ as the left ideal generated by
\begin{equation}\label{eq:defepsilon}
  \epsilon_I=\eta_{\alpha(I)}\xi_{\alpha_{(\bar I^\sim)}}
\end{equation}
where $J^\sim$ denotes the conjugate of a composition $J$.

For a module $M$ over $H_n(0)$, let us say that $M$ is a combinatorial
module if there exists a basis $m_j$ of M such that $\eta_im_j$ is either
$0$ or some $m_k$ (this generalizes the notion of permutation
representation of a group).

Projective $H_n(0)$ modules are combinatorial. The relevant bases are
subsets of a basis of $H_n(0)$ which can also be found in \cite{No}.
Here we will denote it by $g_\sigma$. We set
$g_{\alpha(I)}=\epsilon_I$, and if $\sigma=\tau\alpha(I)$ with
$\ell(\sigma)=\ell(\tau)+\ell(\alpha(I))$ and $\sigma\in D_I$,
\begin{equation}\label{eq:defg}
  g_\sigma= \eta_\tau \epsilon_I \,.
\end{equation}

It is important to mention that the generators $\epsilon_I$ are not
idempotents. The corresponding orthogonal idempotents are denoted by
$e_I$. One way to compute them is to express the identity of $H_n(0)$
in the basis $g_\sigma$. If
\begin{equation}
  1=\sum_{\sigma\in\SG_n}a_\sigma g_\sigma
\end{equation}
then
\begin{equation}
  e_I = \sum_{\sigma\in [\alpha(I),\omega(I)]}
           a_\sigma g_\sigma\,.
\end{equation}

\subsection{Quantum quasi-symmetric functions and quantum shuffles}

It is known that $H_n(0)$-modules are endowed with a natural
filtration, which can be taken into account in the description of the
composition factors of the induced modules
\begin{equation}
  \SS_I\,\widehat{\otimes}\,\SS_J
    =  \SS_I\otimes\SS_J\uparrow_{H_n(0)\otimes H_m(0)}^{H_{m+n}(0)}\,.
\end{equation}
The multiplicity $c_{IJ}^K$ of $\SS_K$ as a composition factor of this
module is equal to the coefficient of $F_K$ in the product $F_IF_J$.
The rule to evaluate this product is as follows: take any permutation
$u$ of $1,\ldots,n$ with descent composition $C(u)=I$ and any
permutation $v$ of $n+1,\ldots, n+m$ such that $C(v)=J$. Then the
shuffle of the two words $u$ and $v$ is a sum of permutations of
$\{1,\ldots, n+m\}$
\begin{equation}
  u\shuffle v  =  \sum_{w\in\SG_{m+n}} c_w w
\end{equation}
and the product is given by
\begin{equation}\label{eq:prodF}
  F_I F_J  =  \sum_{w\in\SG_{m+n}} c_w F_{C(w)} \,.
\end{equation}
There exists a $q$-analogue of the shuffle product, 
which is defined by
\begin{equation}\label{eq:qshuffle}
  \text{if $u=au'$ and $v=bv'$ with $a,b\in A$,\quad then \quad}
  u\shuffle_q v = a(u' \shuffle_q v) +  q^{|u|} b(u\shuffle_q v')
\end{equation}
where $|u|$ is the length of $u$. It can be shown that this operation
is associative, and that when $q$ is not a root of unity, the
$q$-shuffle algebra is isomorphic to the concatenation algebra, which
corresponds to the case $q=0$ \cite{NCSF3}.

The induced representation $\SS_I\,\widehat{\otimes}\,\SS_J$ is
generated by a single vector $u$. There is a filtration of this module
whose $k$-th slice $M_k$ is spanned by the elements $T_\sigma u$ for
permutations $\sigma$ of length $k$. Now, if one computes the product
$F_I F_J$ by using the $q$-shuffle instead of the ordinary one in
formula (\ref{eq:prodF}), the coefficient of $q^kF_H$ in the result is
the multiplicity of $\SS_H$ at level $k$ of the filtration. The
algebra $QSym_q$ of quantum quasi-symmetric functions is defined
accordingly as the algebra with generators $F_I$ and multiplication
rule
\begin{equation}\label{eq:prodqF}
  F_I F_J = \sum_{w} \pair{w}{u\shuffle_q v}F_{C(w)}
\end{equation}
for permutations $u$ and $v$ as above, $\pair{w}{u\shuffle_q v}$ being
the coefficient of $w$ in $u\shuffle_q v$.

All the usual bases of $QSym$, in particular $(M_I)$, are defined in
$QSym_q$ by the same expressions in terms of the $F_I$ as in the
classical case.

For generic values of $q$, $QSym_q$ is freely generated by the
one-part quasi-ribbons $F_n$, or as well by the power-sums $M_n$, or
any sequence corresponding to a free set of generators of the algebra
of symmetric functions in the classical case. This means that if we
define for a composition $I=(i_1,\ldots,i_r)$
\begin{equation}
  F^I = F_{i_1}F_{i_2}\cdots F_{i_r}\quad\text{and}\quad
  M^I = M_{i_1}M_{i_2}\cdots M_{i_r}\ \in \ QSym_q
\end{equation}
then the $F^I$ (resp. the $M^I$) form a basis of $QSym_q$. This
is clearly not true for $q=1$, for in this case these elements
are symmetric functions.

Thus, for generic $q$, $QSym_q$ is isomorphic to the algebra of
noncommutative symmetric functions. Actually, it can be obtained by
specializing the formal variables of the polynomial realization
$\Sym(A)$ of $\Sym$ (see \cite{NCSF1}, Sec.~7.3) to the generators of
the (infinite dimensional) quantum affine space
$\C_q[X]=\C_q[x_1,x_2,\ldots]$, the associative algebra generated by
an infinite sequence of elements $x_i$ subject to the $q$-commutation
relations
\begin{equation} \label{q-commuting}
  \text{for $j>i$,}\quad x_j x_i = q x_i x_j \,.
\end{equation}
More precisely, let $\mathbf{Sym}(X)$ be the subalgebra of $\C_q[X]$
generated by the specialization $a_i \rightarrow x_i$ of the
noncommutative symmetric functions. Then, $\mathbf{Sym}(X)$ is
isomorphic as an algebra to $QSym_q$, the correspondence being given by
\begin{equation}
  M_I \longleftrightarrow \overline M_I
  =  \sum_{j_1<\cdots<j_r}
     x_{j_1}^{i_1}\cdots x_{j_r}^{i_r}\,.
\end{equation}
That is, if ones defines
\begin{equation}
  \overline F_I=\sum_{J\succeq I}\overline M_J\ ,
\end{equation}
one has for $u$ a permutation of $1,\ldots,n$ and $v$ a permutation of
$n+1,\dots,n+m$
\begin{equation}\label{eq:qprodF}
  \overline F_{C(u)} \overline F_{C(v)}
  =   \sum_w \pair{w}{u\shuffle_q v} \overline F_{C(w)}\,.
\end{equation}
Thus, $QSym_q$ provides a kind of unification of both generalizations
$QSym$ and $\Sym$ of $Sym$.

\subsection{Convolution algebras}\label{sec:conv}

Let $H$ be a bialgebra with multiplication $\mu$ and comultiplication
$\Delta$. The convolution product of two endomorphisms $\phi,\psi$ of
$H$ is given by
\begin{equation}\label{eq:conv}
  \phi\star \psi = \mu\circ (\phi\otimes\psi)\circ\Delta \,.
\end{equation}
This is an associative operation, as soon as $\mu$ is associative and
$\Delta$ coassociative. Actually, (\ref{eq:conv}) makes sense, and
is still associative, without
assuming any compatibility between $\mu$ and $\Delta$, and such
expressions will arise in the sequel. When no bialgebra structure is
assumed, we speak of pseudo-convolution.

Interesting examples of convolution algebras are provided by the centralizer
algebras of group actions on tensor spaces. Let $V$ be a
representation of some group $G$. Then, the tensor algebra $T(V)$ is a
representation of $G$, and one can consider its centralizer algebra
\begin{equation}
  H={\rm End}_G T(V) \,.
\end{equation}
It is clearly stable under composition, but also under convolution
since
\begin{align*}
  (\phi\star\psi)(gx)&=\mu\circ(\phi\otimes\psi)\circ\Delta(gx)\\
  &= \mu\circ(\phi\otimes\psi)(g\otimes g)\Delta(x)\\
  &= \mu\circ(g\otimes g)\circ(\phi\otimes\psi)\circ\Delta(x)\\
  &= g (\phi\star\psi)(x)\,.
\end{align*}

When one takes $G=GL(N,\C)$ and $V=\C^N$, $H$ is a homomorphic image
of the direct sum $\C\SG$ of all $\C\SG_n$. By letting
$N\rightarrow\infty$, one obtains a convolution structure on $\C\SG$.
The resulting algebra has been extensively studied by Reutenauer and
his students \cite{Re,MR1,PR}. In the following, we will propose a new
approach leading to a generalization of this algebra.

             \section{Free quasi-symmetric functions}

Our first generalization is obtained by lifting the multiplication
rule (\ref{eq:prodF}) to the free associative algebra, where it
becomes multiplicity free. One arrives in this way to an algebra with
basis labelled by all permutations, which turns out to be isomorphic
to the algebra studied by Malvenuto and Reutenauer in \cite{MR1},
Sec. 3.

\subsection{Free quasi-symmetric functions in a free algebra}

\begin{definition}
The {\em free quasi-ribbon} $\F_\sigma$ labelled by
a permutation  $\sigma\in\SG_n$ is the noncommutative polynomial
\begin{equation}
  \F_\sigma = \sum_{\std(w)=\sigma^{-1}} w \quad \in \Z\free{ A}
\end{equation}
where $\std(w)$ denotes the standardized of the word $w$.
\end{definition}
The hypoplactic version of the Robinson-Schensted correspondence shows
that the commutative image of $\F_\sigma$ is the quasi-symmetric
function $F_I$, where $I=C(\sigma)$. Indeed, the standard ribbon
playing the role of the insertion tableau is equal to $\std(w)^{-1}$,
so that $\F_\sigma$ contains exactly one representative of each
hypoplactic class of shape $I$.

For a word $w=x_1x_2\cdots x_n$ in the letters $1,2,\ldots $ and an
integer $k$, denote by $w[k]$ the shifted word $(x_1+k)(x_2+k)\cdots
(x_n+k)$, e.g., $312[4]=756$. The shifted concatenation of two words
$u$, $v$ is defined by
\begin{equation}
  u\bullet v= u\cdot v[k]
\end{equation}
where $k$ is the length of $u$.

\begin{proposition}\label{MULTF}
Let $\alpha\in\SG_k$ and $\beta\in\SG_l$. Then,
\begin{equation}
  \F_\alpha \F_\beta = \sum_{\sigma\in \alpha\shuffle \beta [k]}\F_\sigma
\end{equation}
Therefore, the free quasi-ribbons span a $\Z$-subalgebra of the free
associative algebra.
\end{proposition}

\Proof
A word $w=a_{i_1}a_{i_2}\ldots a_{i_n}$ can be represented by a
monomial in commuting ``biletters'' $\binom {a_i}{j}$ (which are just
a convenient notation for doubly indexed indeterminates $x_{ij}$). We
identify $w$ with any monomial
$\binom{a_{i_1}}{j_1}\binom{a_{i_2}}{j_2} \cdots \binom{a_{i_n}}{j_n}$
such that $j_1<j_2<\ldots < j_n$, and in particular with the product
$\binom{a_{i_1}}{1}\binom{a_{i_2}}{2} \cdots \binom{a_{i_n}}{n}$,
which we also denote by
\begin{equation}
  \binom{a_{i_1}a_{i_2}\cdots a_{i_n}}{12\cdots n}
  = \binom{w}{{\rm id}}=\binom{w'}{\tau}
\end{equation}
whenever $\tau$ is a permutation such that $w'\tau=w$. Such a
representation if of course not unique. Then, $\sigma=\std(w)^{-1}$ is
the unique permutation of minimal length such that $\binom{w}{{\rm
id}}=\binom{w^+}{\sigma}$, where $w^+$ denotes the non-decreasing
rearrangement of $w$. The correspondence $w\leftrightarrow
\binom{w^+}{\sigma}$ is a bijection between words and pairs
$(u,\alpha)$ where $u$ is a nondecreasing word and $\alpha$ is a
permutation of the same length such that $\alpha_i<\alpha_{i+1}$ when
$u_i=u_{i+1}$. In this case, we say that $u$ is $\alpha$-compatible,
and we write $u\uparrow\alpha$. The concatenation product corresponds
to an operation $\circ$ on biwords, given by the rule
\begin{equation}
  \binom{u}{\alpha}\circ\binom{v}{\beta}
  = \binom{uv}{\alpha\bullet\beta}\,.
\end{equation}
Now, if we write
\begin{equation}
  \F_\alpha = \sum_{u\uparrow\alpha}\binom{u}{\alpha}\,,\qquad
  \F_\beta = \sum_{v\uparrow\beta}\binom{v}{\beta}\,,
\end{equation}
where $u$ and $v$ run over nondecreasing words of respective
lengths $k$ and $l$, we see that
\begin{equation}
  \F_\alpha \F_\beta
  =
  \sum_{u\uparrow\alpha,v\uparrow\beta}
    \binom{uv}{\alpha\cdot\beta[k]}
  =
  \sum_{\sigma\in\alpha\shuffle\beta[k]}
    \sum_{w\uparrow\sigma}
      \binom{w}{\sigma}\,,
\end{equation}
whence the proposition. \qed

\begin{definition}
The subalgebra of $\C \free{A}$
\begin{equation}
  \FQSym = \bigoplus_{n\ge 0}\bigoplus_{\sigma\in\SG_n}  \C\,\F_\sigma
\end{equation}
is called the algebra of {\em free quasi-symmetric functions\/}.
\end{definition}

It will be convenient to define a scalar product on $\FQSym$
by setting
\begin{equation}
  \pairing{\F_\sigma}{\F_\tau}=\delta_{\sigma^{-1},\tau}
\end{equation}
and to introduce the  notation
\begin{equation}
  \G_\sigma=\F_{\sigma^{-1}}
\end{equation}
for the adjoint basis of $(\F_\sigma)$.

Since the convolution of permutations is related to
the shifted shuffle by 
\begin{equation}
  (\alpha^{\vee\ }\star\beta^{\vee\ })^{\vee\ } 
      = \alpha \shuffle \beta[k]
\end{equation}
where $f\to f^{\vee\ }$ is the linear involution defined on
permutations by $\alpha\to\alpha^{\vee\ }=\alpha^{-1}$,  
we see that $\FQSym$ is isomorphic to the convolution algebra
of permutations of \cite{MR1}. The interesting point is that
the natural map $\sigma\mapsto F_{C(\sigma)}$ from this
algebra to $QSym$ becomes simply the commutative image $a_i\mapsto x_i$.

The quasi-symmetric generating function of a set of permutations in
the sense of \cite{Ge} can now be regarded as the commutative image of
an element of $\FQSym$. We shall see that in certain special cases,
such as linear extensions of posets, the free quasi-symmetric function
can be more interesting ({\it cf}. Section \ref{Sec:Hn0-Poset}).

Another  property of $\FQSym$ is that it contains
a subalgebra with a distinguished basis  labelled by standard
Young tableaux, which maps to ordinary Schur functions under
abelianization, and to which the Littlewood-Richardson
rule can be lifted to a multiplicity free formula (see Proposition \ref{LRS}).                                
The kind of argument used to establish this formula can also be used to prove
Proposition \ref{MULTF}.  
Recall that we  denote by $w\mapsto (\HQ(w),\HR(w))$ the hypoplactic
Robinson-Schensted correspondence. An alternative
definition of $\F_\sigma$ is
\begin{equation}
  \F_\sigma=\sum_{\HR(w)=\sigma} w
\end{equation}
and Proposition \ref{MULTF} can be derived exactly in the same way as
Proposition \ref{LRS},
from the fact that the hypoplactic congruence is compatible to
restriction to intervals (see \cite{Loth}).

\subsection{Duality}
 
One can define on $\FQSym$ a bialgebra structure imitated from the
case of ordinary quasi-symmetric functions.
 
Let $A'$ and $A''$ be two mutually commuting ordered alphabets.
Identifying $F\otimes G$ with $F(A')G(A'')$, we set
$\Delta(F)=F(A'\oplus A'')$, where $\oplus$ denotes the ordered sum.
Clearly, this is an algebra homomorphism.

\begin{proposition}
  $\FQSym$ is a bialgebra for $\Delta$, and on the
  basis $\F_\sigma$, the comultiplication is given by
  \begin{equation}
    \Delta\F_\sigma=\sum_{u\cdot v=\sigma}\F_{\std(u)}\otimes\F_{\std(v)}\,.
  \end{equation}
  where $u\cdot v$ denotes the concatenation of $u$ and $v$.
\end{proposition}

\Proof 
Like Proposition \ref{MULTF}, this formula is easily obtained
in the biword notation. Indeed, $\F_\sigma(A'\oplus A'')$ is the
image of the element
\begin{equation}\label{eq:comp}
  \sum_{w\uparrow\sigma}{\binom{w}{\sigma}}=
  \sum_{w=w'w''\uparrow\sigma}{\binom{w'w''}{\sigma}}
\end{equation}
of the free algebra $\C\free{A'\cup A''}$ under the map
\begin{equation}
  \pi : \C\free{A'\cup A''}
  \rightarrow
  \C\free{A',A''}
  \simeq
  \C\free{A}\otimes\C\free{A}.
\end{equation}
The sum runs over all nondecreasing
$\sigma$-compatible words $w$, which are necessarily of the form
$w'w''$ with $w'\in {A'}^*$ and $w''\in {A''}^*$, since $A'<A''$. Let
$k=|w'|$ and $l=|w''|$. As a word, $\sigma$ can be factorized as
$\sigma=uv$, where $|u|=k$ and $|v|=l$. Let $\sigma'=\std(u)$ and
$\sigma''=\std(v)$. Since $A'$ and $A''$ are disjoint, $w'$ and $w''$
have to be respectively $\sigma'$ and $\sigma''$ compatible, and
actually, the sum (\ref{eq:comp}) runs exactly over all such words.
Since
\begin{equation}
  \pi \binom{w'w''}{\sigma}=\binom{w'}{\sigma'}\circ \binom{w''}{\sigma''}
\end{equation}
the image under $\pi$ of this sum factorizes into
\begin{equation}
  \sum_{uv=\sigma}\sum_{w'\uparrow\std(u)}\binom{w'}{\std(u)}
                  \sum_{w''\uparrow\std(v)}\binom{w''}{\std(v)}
\end{equation}
whence the proposition. \qed
\bigskip
 
\begin{corollary} $\FQSym$ is a self-dual bialgebra. That is,
  for all $F,G,H\in\FQSym$,
  \begin{equation}  
    \pairing{F\otimes G}{\Delta H}=\pairing{FG}{H}\,.
  \end{equation}    
\end{corollary}
 
\Proof Denote by $\odot$ the multiplication adjoint to $\Delta$,
that is, such that
\begin{equation}
  \pairing{F\otimes G}{\Delta H}=\pairing{F\odot G}{H}\,,
\end{equation}
and consider the structure constants
\begin{equation}
  \G_\alpha\odot \G_\beta =\sum_\gamma g_{\alpha\beta}^\gamma G_\gamma\,.
\end{equation}
Then, $ g_{\alpha\beta}^\gamma =1$ if $\alpha=\std(u)$ and
$\beta=\std(v)$ for some factorization $\gamma=uv$, and
$g_{\alpha\beta}^\gamma =0$ otherwise. Therefore, these structure
constants coincide with those of the convolution product on
permutations:
\begin{equation}
  \alpha\star\beta = \sum_\gamma   g_{\alpha\beta}^\gamma \gamma\,.
\end{equation}
\bigskip

We have therefore interpreted the two multiplications and comultiplications
of \cite{MR1} as operations on labels of two different bases of the same
subalgebra of the free associative algebra.

\subsection{Algebraic structure}
 
We can now apply to $\FQSym$ the results of Poirier and Reutenauer
\cite{PR} and we see that $\FQSym$ is freely generated by the
$\G_\sigma$, where $\sigma$ runs over {\em connected permutations}
(see \cite{Comtet}), \ie  permutations such that $\sigma([1,k])\not =
[1,k]$ for all intervals $[1,k]\subseteq [1,n-1]$. Actually, this
result holds for a one-parameter family of algebras, and we shall now
reprove it in this context.

We denote by $\conn$ the set of connected permutations, and by
$c_n=|\conn_n|$ the number of such permutations in $\SG_n$. For later
reference, we recall that the generating series of $c_n$ is
\begin{align*}
  \sum_{n\ge 1} c_n t^n
  =& 1 - \left(\sum_{n\ge 0} n! t^n\right)^{-1}\\
  =& t+{t}^{2}+3\,{t}^{3}+13\,{t}^{4}+71\,{t}^{5}+461\,{t}^{6}
     + 3447\,{t}^{7 }+29093\,{t}^{8}\\
  &+273343\,{t}^{9}+2829325\,{t}^{10}+31998903\,{t}^{11}
     + 392743957\,{t}^{12} + O(t^{13})\,.
\end{align*}

For $\alpha\in\SG_k$ and $\beta\in\SG_l$, recall that
$\alpha\bullet\beta= \alpha\cdot\beta[k]$ is the shifted concatenation
of $\alpha$ and $\beta$. Any permutation $\sigma\in\SG_n$ has a unique
maximal factorization $\sigma=\sigma_1\bullet\cdots\bullet\sigma_r$
into connected permutations. Then, the elements
\begin{equation}
  \G^\sigma=\G_{\sigma_1}\cdots\G_{\sigma_r}
\end{equation}
and
\begin{equation}
  \F^\sigma=\F_{\sigma_1}\cdots\F_{\sigma_r} 
\end{equation}
form two bases of $\FQSym$. Since $(\alpha^{-1}\bullet\beta^{-1})^{-1}
=\alpha\bullet\beta$, we have $\F^\sigma=\G^{\sigma^{-1}}$, and the
multiplication of $\FQSym$ is given in both bases by the same formula:
$\G^\alpha\G^\beta=\G^{\alpha\bullet\beta}$ and
$\F^\alpha\F^\beta=\F^{\alpha\bullet\beta}$.

The operations on permutations $\alpha\bullet\beta$ and
$\alpha\shuffle\beta[k]$ describing the multiplication in the bases
$\F^\sigma$ and $\F_\sigma$ are the cases $q=0$ and $q=1$ of the
shifted $q$-shuffle $\alpha\shuffle_q\beta[k]$. This suggests the
consideration of a $q$-deformed algebra $\FQSym_q$, defined as the
(abstract) algebra with generators $\ff_\sigma$ and relations
$\ff_\alpha\ff_\beta=\ff_{\alpha\shuffle_q\beta[k]}$ (where linearity
of the symbol $\ff$ with respects to subscripts is understood). As
above, let
$\ff^\sigma=\ff^{\sigma_1}\cdots\ff^{\sigma_r}=\ff_\sigma+O(q)$. For
each $n$, the $n!\times n!$ matrix expressing the elements
$\ff^\sigma$ on the basis $\F_\sigma$ is of the form $I+O(q)$, and is
therefore invertible over $\C[[q]]$. Moreover, it is unitriangular
with respect to the lexicographic order on permutations, so that it is
actually invertible over $\C[q]$. This proves that the algebras
$\FQSym_q$ are actually isomorphic to each other for all values of
$q$. For $q\not=0$, the isomorphism $\FQSym\rightarrow \FQSym_q$ is
realized by $\F_\sigma\mapsto q^{l(\sigma)}\ff_\sigma$, and for $q=0$,
by $\F_\sigma\mapsto \ff^\sigma$.

\subsection{Primitive elements}
 
Let $\L$ be the primitive Lie algebra of $\FQSym$.
Since $\Delta$ is not cocommutative, $\FQSym$ cannot be the universal
enveloping algebra of $\L$. Let $l_n = \dim \L_n$.

Let us recall that $\G^\sigma=\G_{\sigma_1}\cdots\G_{\sigma_r}$ where
$\sigma=\sigma_1\bullet\cdots\bullet\sigma_r$ is the unique maximal
factorization of $\sigma\in\SG_n$ into connected permutations.
\begin{proposition} 
Let $\V_\sigma$ be the adjoint basis of $\G^\sigma$.
Then, the family $(\V_\alpha)_{\alpha\in\conn}$
is a basis of $\L$. In particular, we have $l_n=c_n$.
\end{proposition}

\Proof If $\alpha$ is connected, then
\begin{equation*}
  \begin{split}
    \Delta\V_\alpha & = \sum_{\sigma,\tau}
    \pairing{\Delta\V_\alpha}{\G^\sigma\otimes\G^\tau}
    \V_\sigma\otimes\V_\tau \\
                    & = \sum_{\sigma,\tau}
    \pairing{\V_\alpha}{\G^{\sigma\bullet\tau}}\V_\sigma\otimes\V_\tau
    =\V_\alpha\otimes 1 + 1\otimes \V_\alpha
  \end{split}
\end{equation*}
since the only possible factorization of $\alpha$ is
$\alpha=\emptyset\bullet\alpha=\alpha\bullet\emptyset$, where
$\emptyset$ denotes the empty word.

Conversely, let $Z=\sum_\alpha c_\alpha \V_\alpha$ be a primitive
element. If $\alpha$ is not connected, let $\alpha=\sigma\bullet\tau$
be a non-trivial factorization. Then,
\begin{equation}
\pairing{\Delta Z}{\G^\sigma\otimes \G^\tau}
  =  \pairing{Z }{\G^{\sigma\bullet\tau}}
  =  \pairing{Z}{\G^\alpha}=c_\alpha
\end{equation}
which has to be zero since the left-hand side is the coefficient
of $\V_\sigma\otimes \V_\tau$ in $\Delta Z$. \qed

\begin{example} In degree 3 we have
\begin{eqnarray*}
\mathbf{V}_{{312}}&=&\F_{{312}}-\F_{{213}}\\
\mathbf{V}_{{231}}&=&-\F_{{132}}+\F_{{231}}\\
\mathbf{V}_{{321}}&=&\F_{{123}}-\F_{{132}}-\F_{{213}}+\F_{{321}}
\end{eqnarray*}

and in degree 4
\begin{eqnarray*}
\mathbf{V}_{{4123}}&=&\F_{{4123}}-\F_{{3124}}\\
\mathbf{V}_{{4132}}&=&\F_{{4132}}-\F_{{3124}}+\F_{{2134}}-\F_{{2143}}\\
\mathbf{V}_{{3412}}&=&-\F_{{1423}}+\F_{{1324}}+\F_{{3412}}-\F_{{2314}}\\
\mathbf{V}_{{3142}}&=&\F_{{3142}}-\F_{{2143}}\\
\mathbf{V}_{{4312}}&=&-\F_{{1423}}+\F_{{1324}}+\F_{{4312}}-\F_{{3214}}\\
\mathbf{V}_{{2413}}&=&-\F_{{1423}}+\F_{{1324}}+\F_{{2413}}-\F_{{2314}}\\
\mathbf{V}_{{4213}}&=&\F_{{4213}}-\F_{{3214}}\\
\mathbf{V}_{{2431}}&=&-\F_{{1432}}+\F_{{2431}}\\
\mathbf{V}_{{2341}}&=&-\F_{{1342}}+\F_{{2341}}\\
\mathbf{V}_{{4231}}&=&\F_{{1243}}-\F_{{1342}}-\F_{{3124}}+\F_{{2134}}-\F_{{2143}}+\F_{{4231}}\\
\mathbf{V}_{{3421}}&=&\F_{{1324}}-\F_{{1432}}-\F_{{2314}}+\F_{{3421}}\\
\mathbf{V}_{{3241}}&=&\F_{{1243}}-\F_{{1342}}-\F_{{2143}}+\F_{{3241}}\\
\mathbf{V}_{{4321}}&=&-\F_{{1234}}+\F_{{1243}}+\F_{{1324}}-\F_{{1432}}+\F_{{2134}}-\F_{{2143}}-\F_{{3214}}+\F_{{4321}}
\end{eqnarray*}
\end{example}

The Hilbert series of the universal enveloping algebra $U(\L)$
(the domain of cocommutativity of $\Delta$) is
\begin{equation*}
  \begin{split}
    \prod_{n\ge 1}(1-t^n)^{-c_n} &= 
    1+t+2\,{t}^{2}+5\,{t}^{3}+19\,{t}^{4}+93\,{t}^{5}\\
    &+574\,{t}^{6}+4134\, {t}^{7}+34012\,{t}^{8}+313231\,{t}^{9}+3191402\,{t}^{10}\\
    &+35635044\,{t}^{11}+432812643\,{t}^{12}+O\left ({t}^{13}\right )\,.
  \end{split}
\end{equation*}

\begin{conjecture}
  $\L$ is a free Lie   algebra.
\end{conjecture}

Assuming the conjecture, denote by $d_n$ the number of generators
of degree $n$ of $\L$. Then, using the
$\lambda$-ring notation, since $\sigma_1\circ L =(1-p_1)^{-1}$
(where $\sigma_1=\sum_{n\ge 0}h_n$, $L=\sum_{n\ge_1} \ell_n$, and
$\ell_n=\frac{1}{n}\sum_{d|n}\mu(d)p_d^{n/d}$ are the Lie characters,
or Witt symmetric functions), we have the equivalent plethystic equations
\begin{equation}
  L\left[\sum_{n\ge 1} d_n t^n\right]=\sum_{n\ge 1} c_n t^n\,,
\end{equation}
and
\begin{equation}
  d(t)= \sum_{n\ge 1} d_n t^n
      = 1-\lambda_{-1}\left[\sum_{n\ge 1}c_n t^n\right]
      = 1-\prod_{n\ge 1}(1-t^n)^{c_n}\,.
\end{equation}
Numerical calculation gives for the first terms
\begin{equation*}
  \begin{split}
    d(t) & = t+{t}^{2}+2\,{t}^{3}+10\,{t}^{4}+55\,{t}^{5}+377\,{t}^{6}\\
         & + 2892\,{t}^{7}+25007\,{t}^{8}+239286\,{t}^{9}+2514113\,{t}^{10}
           + 28781748\,{t}^{11}\\
         & + 356825354\,{t}^{12} + O(t^{13})\,.
  \end{split}
\end{equation*}

\bigskip

We shall now give a formula for the projector $\pi:\ \FQSym\rightarrow
\L$ such that
\begin{equation*}
  \pi(\F_\alpha) = \left\{
    \begin{array}{cl}
      0              & \ \hbox{if $\alpha$ is not connected} \\[1mm]
      \mathbf{V}_\alpha & \ \hbox{if $\alpha$ is connected} \ . 
    \end{array}
  \right.
\end{equation*}

\medskip Let $p_n$ denote the projection onto the homogeneous
component $\FQSym_n$ of $\FQSym$, and let $\mu_q:\ 
\F_\alpha\otimes\F_\beta \mapsto \F_{\alpha\shuffle_q \beta[k]}$ be
the multiplication map of $\FQSym_q$. The $q$-convolution of two
graded linear endomorphisms $f,g$ of $\FQSym$ is defined by
\begin{equation}
f \odot_q g = \mu_q\circ(f\otimes g)\circ\Delta\, .
\end{equation}
For $q=1$, this reduces to ordinary convolution,
otherwise, it is an example of pseudo-convolution
as defined in \ref{sec:conv}. We shall be interested
in the case $q=0$. For a composition $I=(i_1,\ldots,i_m)$, let
\begin{equation}
p_I= p_{i_1}\odot_0\cdots\odot_0 p_{i_m} \,.
\end{equation}

\begin{lemma}
  The $p_I$ are mutually commuting projectors. More precisely we have
  \begin{equation*}
    p_I\circ p_J= \left\{
      \begin{array}{cl}
        0  & \ \hbox{if} \ |I|\not= |J|. \\[1mm] 
        p_{I\vee J} & \ \hbox{otherwise} \ .
      \end{array}
    \right.
  \end{equation*}
  where $I\vee J$ is the composition with descent set
  $\Des(I)\cup\Des(J)$.
\end{lemma}
\Proof The result is clear when $|I|\not= |J|$. Otherwise, we suppose
$|I|= |J|=n$ and proceed by induction on $d=\min (l(I\vee
J)-l(I),l(I\vee J)-l(J))$. If $d=0$ it is easy to check that $p_I\circ
p_{I\vee J}=p_{I\vee J}\circ p_I=p_{I\vee J}$ otherwise, the induction
step is a consequence of the standardization inertia $\std(\sigma
\bullet \tau)=\std(\sigma) \bullet \std(\tau)$ \cqfd
\medskip

Before stating the main proposition we need some notation: For a word
of length $n$, $w=a_1a_2\cdots a_n$ and $S=\{s_1,s_2\cdots
s_k\}\subset [1..n]$ a subset in increasing order, we denote the
corresponding subword by $w|_S=a_{s_1}a_{s_2}\cdots a_{s_k}$. Let
$I=(i_1,i_2,\cdots i_m)$ be a composition of weight $n$. The
factorization-standardization operator $\sfact_I$ is defined by
\begin{equation*}
  \sfact_I(w)=
  \left\{
    \begin{array}{cl}
      \std(w|_{[1..i_1]})\tensor \std(w|_{[i_1+1..i_1+i_2]})
      \tensor\cdots \std(w|_{[n-i_m+1..n]}) &
      \ \hbox{if} \ |w|=n, \\[1mm]
      0  & \ \hbox{otherwise}
    \end{array}
  \right.
\end{equation*}
For example
$\sfact_{(2,3)}(53412)=\std(53)\tensor\std(412)=21\tensor312$. 
We can now state:
 
\begin{proposition} 
  \begin{itemize}
  \item[(i)] The operator
    \begin{equation}
      \displaystyle\pi=\sum_{|I|\ge 1} (-1)^{l(I)-1} p_I
\end{equation}
is the projector onto the primitive Lie algebra with the span of
    $(\F_\alpha)_{\alpha \notin\conn}$ as kernel.
  \item[(ii)] Moreover, one has $\mathbf{V}_\alpha = \pi(\F_\alpha)$ for
    $\alpha$ connected.
  \end{itemize}
\end{proposition}
 
\Proof                                                                    
The m-fold shifted concatenation ${\rm sconc}^{(m)}$ is defined in the obvious
way. Then, for $l(I)=m$, 
\begin{equation*}
  p_I(\F_\alpha)=
  \left\{
    \begin{array}{cl}
      \F_{{\rm sconc}^{(m)}\,\circ\,\sfact_I(\alpha)}
        & \text{if $\alpha\in\SG_n$\,,}
      \\[1mm]
      0 & \text{otherwise.}
    \end{array}
\right.
\end{equation*}                                                                        

We first prove that, if $l(I_0)=2$, one has
$\pi\circ p_{I_0}=0$. For $i=0,1$,
let $${\cal I}_i=\{|I|=n\ |\ \delta(I_0\prec I)=i\}\,,$$ it is easy
to check that $\sharp ({\cal I}_i)=2^{n-2}$ and that
$I\rightarrow I\vee I_0$ induces a bijection
${\cal I}_0\rightarrow {\cal I}_1$. Hence
\begin{equation*}
  \begin{split}
    \pi\circ p_{I_0} &=\sum_{l(I)\geq 1}(-1)^{l(I)-1}p_I\circ p_{I_0}
      = \sum_{|I|=n}(-1)^{l(I)-1}p_{I\vee I_0} \\
    & = \sum_{I\in {\cal I}_0}(-1)^{l(I)-1} p_{I\vee I_0} +
        \sum_{I\in {\cal I}_1}(-1)^{l(I)-1} p_{I\vee I_0} \\
    & = \sum_{I\in {\cal I}_1}(-1)^{l(I)} p_I +
        \sum_{I\in {\cal I}_1}(-1)^{l(I)-1} p_I = 0
  \end{split}
\end{equation*}
                                                                        
If $\alpha\notin \conn_n$, then for some composition $I_0$ of $n$ of
length 2, we have $p_{I_0}(\F_\alpha)=\F_\alpha$. Hence $\F_\alpha\in
\ker(\pi)$. Now, if $\alpha\in \conn_n$, the construction of $\pi$
shows that
\begin{equation}\label{triang}
  \pi(\F_\alpha)=\F_\alpha +
  \sum_{\beta\notin \conn_n}c_\beta \F_\beta
\end{equation}
and then $\pi^2(\F_\alpha)=\pi(\F_\alpha)$. This finishes to prove
that $\pi$ is a projector and from (\ref{triang}), we get that the
generating series of $Im(\pi)$ is exactly $\sum_n c_nt^n$.

The comultiplication on $\FQSym$ can be rewritten as
\begin{equation}
  \Delta=Id\tensor 1+ 1\tensor Id + \sum_{l(I)=2} \sfact_I
\end{equation}
so, to get $\operatorname{Im}(\pi)\subset \L$, it suffices to prove
\begin{equation}
  \left(\sum_{l(I)=2} \sfact_I\right)\circ \pi=0.
\end{equation}

But, from the construction of $\sfact_I$, one has
$\sfact_I=\sfact_I\circ p_I$. Now, if $l(I)=2$, we get
\begin{equation}
  \sfact_I\circ \pi=\sfact_I\circ p_I\circ \pi
  = \sfact_I\circ (\pi\circ p_I)=0
\end{equation}
which proves that $Im(\pi)\subset \L$, the equality
of these two spaces follows from the fact that the generating series are
equal.
 
Equation (\ref{triang}) says that $(\pi(F_\alpha))_{\alpha \in\conn }$
is the unique basis of $\L$ such that
\begin{equation}
  \pi(\F_\alpha)=\F_\alpha + \sum_{\beta\notin \conn_n}c_\beta
\F_\beta\,.
\end{equation}
Since $\V_{\alpha}$ also have this property,
$\V_{\alpha}=\pi(\F_\alpha)$.  \cqfd

\subsection{Free symmetric functions and the Littlewood-Richardson rule}

\begin{definition}
  Let $t$ be a standard tableau of shape $\lambda$. The {\em free
    Schur function} labelled by $t$ is
  \begin{equation}
    \SS_t = \sum_{P(\sigma)=t} \F_\sigma=\sum_{Q(w)=t }w \,,
  \end{equation}
  where $w\mapsto (P(w),Q(w))$ is the usual Robinson-Schensted map.
\end{definition}

As pointed out in \cite{Loth}, Sch\"utzenberger's version of the
Littlewood-Richardson rule is equivalent to the following statement,
which shows that the free Schur functions span a subalgebra of
$\FQSym$. We will call it the algebra of {\em free symmetric
  functions} and denote it by $\FSym$. It provides a realization of
the algebra of tableaux introduced by Poirier and Reutenauer \cite{PR}
as a subalgebra of the free associative algebra. A representation
theoretical interpretation will be given in the sequel.

\begin{proposition}[LRS rule]\label{LRS}
  Let $t'$, $t''$ be standard tableaux, and let $k$ be the number of
  cells of $t'$.  Then,
  \begin{equation}
    \SS_{t'}\SS_{t''} = \sum_{t\in \Sh(t',t'')} \SS_t 
  \end{equation}
  where $\Sh(t',t'')$ is the set of standard tableaux in the shuffle
  of $t'$ (regarded as a word via its row reading) with the plactic
  class of $t''[k]$.
\end{proposition}

\bigskip

\Proof This follows from Proposition \ref{MULTF}, and the fact
that the plactic congruence is compatible with restriction to intervals.
Indeed, denote by $\equiv$ the plactic congruence on the free
algebra $\Z\free{A}$, for some ordered alphabet $A=\{a_1<a_2<\cdots < a_n\}$.
For a word $w\in A^*$ and an interval $I=[a_i,a_j]$ of $A$, denote
by $w|_I$ the word obtained by erasing in $w$ the letters not in $I$.
Then, since the plactic relations $xzy\equiv zxy\ (z\le y<z$
and $yxz \equiv yzx \ (x<y\le z)$ reduce to equalities after  erasing 
$x$ or $z$,  $w\equiv w' \Rightarrow w|_I \equiv w'|_I$.
From this, we see that
\begin{equation}
  \SS_{t'}\SS_{t''}
  = \sum_{P(\sigma')=t, P(\sigma'')=t''}  \F_{\sigma'}\F_{\sigma''} 
  = \sum_{t\in\Sh(t',t'')}\sum_{P(w)=t}  \F_\sigma\,,
\end{equation}
since the set of permutations 
$\{\sigma'\shuffle \sigma''[k]\ |P(\sigma')=t,P(\sigma'')=t''\}$
is, by the above remark, a union of plactic classes. Each class contains
a unique tableau, and since the restriction of a tableau to an initial
segment of the alphabet has to be a tableau, such a tableau can
appear only in the shuffles $t'\shuffle \sigma''[k]$.
\qed

\bigskip
The original Littlewood-Richardson rule, as well as its plactic
version, are immediate corollaries of Proposition \ref{LRS} (see \cite{Loth}).

\begin{example}
  The smallest interesting example occurs for the shape $(2,1)$, e.g.,
  with
  $$
  t'=t''=\ \raisebox{-0.3\height}{\young{3\cr 1&2\cr}}\ ,
  $$
  the product $\SS_{t'}\SS_{t''}$ is equal to $\sum_t \SS_t$
  where $t$ ranges over the following tableaux:
  $$
  \young{3&6\cr 1&2&4&5\cr}\qquad
  \young{3&4&6\cr 1&2&5\cr}\qquad
  \young{6\cr 3\cr 1&2&4&5\cr}\qquad
  \young{4\cr 3&6\cr 1&2&5\cr}
  $$
  $$
  \young{6\cr 3&4\cr 1&2&5\cr}\qquad
  \young{4&6\cr 3&5\cr 1&2\cr}\qquad
  \young{6\cr 4\cr 3\cr 1&2&5\cr}\qquad
  \young{6\cr 4\cr 3&5\cr 1&2\cr}
  $$
\end{example}

The scalar product of two
free Schur functions is equal to 1 whenever the corresponding tableaux
have the same shape, and to 0 otherwise. Indeed,
\begin{equation}
  \pairing{\SS_{t'}}{\SS_{t''}}
  = \sum_{P(\sigma)=t',Q(\sigma)=t''} \pairing{\F_\sigma}{\G_\sigma}
  = 1
\end{equation}
since a permutation is uniquely determined by its $P$ and $Q$ symbols.
\bigskip

Note that the algebra of noncommutative symmetric functions
$\Sym(A)$ is a subalgebra of $\FSym$, since
\begin{equation}
  R_I(A) = \sum_{\Rec(t)=\Des(I)}  \SS_t\,,
\end{equation}
where $\Rec(t)$ denotes the recoil (or descent) set of the tableau $t$.

\subsection{An example: the $QS$-distribution on symmetric groups}
 
The definitions of this section are well illustrated by a certain
probability distribution on symmetric groups investigated by Stanley
in \cite{St}. Let $x=(x_i)_{i\ge 1}$ be a probability distribution on
our infinite alphabet $A=\{a_1,a_2,\ldots\}$, that is,
$\prob(a_i)=x_i$, $x_i\ge 0$, and $\sum x_i=1$. From this, one defines
a probability distribution $QS(x)$ on each symmetric group $\SG_n$ by
the formula
\begin{equation}
  \prob(\sigma)=\G_\sigma (x) \ .
\end{equation}
That this is actually a probability distribution follows from the
identity $\G_1^n=\sum_{\sigma\in\SG_n}\G_\sigma$. Then, Theorem 2.1 of
\cite{St} states that $\prob(\sigma)=F_{C(\sigma^{-1})}(x)$, which
follows from the equalities $\G_\sigma=\F_{\sigma^{-1}}$ and
$\F_\sigma(x)=F_{C(\sigma)}(x)$.
 
Next, Stanley introduces the operator
\begin{equation}
  \Gamma_n(x)=\sum_{\sigma\in\SG_n}\prob(\sigma)\sigma \ \in\ \C\SG_n\ .
\end{equation}
Actually, $\Gamma_n(x)$ is in the descent algebra $\Sigma_n$,
and the corresponding noncommutative symmetric function is $S_n(xA)$.
Therefore \cite{NCSF2}, the eigenvalues of $\Gamma_n(x)$ are
the $p_\lambda(x)$, with multiplicities $n!/z_\lambda$. Also,
the convolution formula $\Gamma_n(x)\Gamma_n(y)=\Gamma_n(xy)$ amounts to
the identity $S_n(xA)*S_n(yA)=S_n(xyA)$ of \cite{NCSF2}.
 
Another result of \cite{St} is that the probability $M_n(k)$ that a random
permutation (chosen from the $QS$-distribution) has $k$ inversions
is equal to the probability that it has major index $k$ (Theorem 3.2).
This is equivalent to the identity
\begin{equation}
  \sum_{\sigma\in\SG_n} q^{l(\sigma)}\G_\sigma(x)
  = \sum_{\sigma\in\SG_n} q^{\maj(\sigma)}\G_\sigma(x) .
\end{equation}
The right-hand side can be rewritten as
\begin{multline*}
  \sum_{|I|=n}q^{\maj(I)}r_I(x)
  = \sum_{I,J}q^{\maj(I)}\pairing{r_I}{r_J}F_J(x) \\
  = \sum_J\left(\sum_{C(\sigma^{-1})=J}q^{\maj(\sigma)}\right) F_J(x)
  = \sum_J\left(\sum_{C(\sigma^{-1})=J}q^{l(\sigma)}\right)F_J(x)
  = \sum_{\sigma\in\SG_n} q^{l(\sigma)}\G_\sigma(x) 
\end{multline*}
since $l(\sigma)$ and $\maj(\sigma^{-1})$ have the same
distribution on a descent class ({\it cf. } \cite{St1,FS}).
 
Finally, we note that the specialization $\SS_t(x)$ of a free Schur
function is the probability that a $QS$-random permutation has $t$ as
insertion tableau, and that $R_I(x)$ is the probability that a random
permutation has shape $I$ (Theorems 3.4 and 3.6 of \cite{St}).

\subsection{Quantum quasi-symmetric functions again}

Recall that we denote by $\C_q[X]$ the algebra of quantum polynomials,
generated by letters $x_i$ subject to the relations
$x_jx_i=qx_ix_j$ when $j>i$. The following proposition clarifies
the constructions of \cite{TU}.

\begin{proposition}\label{prop:qsf}
The natural homomorphism $\varphi_q:\ a_i\mapsto x_i$ from $\C\free{A}$ to
the algebra of quantum polynomials
$\C_q[X]$ maps $\F_\sigma$ to the quantum quasi-symmetric
function $q^{\ell(\sigma)}F_{C(\sigma)}$.
\end{proposition}

\Proof For any word $w$, one has
$$
\varphi_q(w)=q^{\ell(\sigma)}\varphi_q(w^+)
$$
where $\sigma=\std(w)^{-1}$ and $w^+$ is the nondecreasing
rearrangement of $w$. \qed

\bigskip 
Therefore, $QSym_q$ is a quotient of $\FQSym$.
The multiplication formula (\ref{eq:prodqF}) appears now as an
immediate consequence of Proposition \ref{MULTF}.
The $q$-generating function $\Gamma_q(P)$ of a poset, introduced
in \cite{TU} to derive  (\ref{eq:prodqF}), 
can also be regarded as the image under $\varphi_q$ of 
of a free generating function $\bGamma(P)$ described in the forthcoming
section.
Most formulas of \cite{TU} are easy consequences of Proposition \ref{prop:qsf}.
For example, formula (38) of \cite{TU}, which can be stated as
\begin{equation}
  \varphi_q(R_I(A))=\sum_J c_{IJ}(q)\overline{F}_J
\end{equation}
where 
\begin{equation}
  c_{IJ}(q)=\sum_{C(\sigma)=I,\ C(\sigma^{-1})=J}q^{\ell(\sigma)}
\end{equation}
follows from the expression
\begin{equation}
  R_I(A)=\sum_{C(\sigma)=I}\G_\sigma\,.
\end{equation}
\bigskip

\subsection{Posets, $P$-partitions, and the like}
\label{Sec:Hn0-Poset}

Here, by a {\em poset}, we  mean any partial order $P$ on the
set $[n]=\{1,2,\ldots,n\}$. We write $<_P$ for the order of $P$
and $<$ for the usual total order on $[n]$.  Stanley \cite{St1} 
defines a $P$-partition as a function $f:\ [n]\rightarrow X$
for some totally ordered set of variables $X$, such that
\begin{equation}
  i <_P j \Rightarrow\ f(i)\le f(j)\quad{\rm and}\quad
  i <_P j\ {\rm and}\ i>j\ \Rightarrow f(i)<f(j)\,.
\end{equation}
In \cite{Ge},  Gessel associates to a poset $P$ a generating
function
\begin{equation}
  \Gamma(P)=\sum_{f\in {\cal A}(P)} f(1)f(2)\cdots f(n)
\end{equation}
where ${\cal A}(P)$ denotes the set of all $P$-partitions.
This generating function turns out to be quasi-symmetric,
actually,
\begin{equation}
  \Gamma(P) = \sum_{\sigma\in L(P)} \Gamma(\sigma)
            = \sum_{\sigma\in L(P)} F_{C(\sigma)}
\end{equation}
where $L(P)$ denotes the set of linear extensions of $P$, which can be
identified with permutations $\sigma\in\SG_n$ such that $i<_P
j\Rightarrow \sigma^{-1}(i)<\sigma^{-1}(j)$. Identifying a
$P$-partition $f$ with the word $w_f=a_{f(1)}a_{f(2)}\cdots a_{f(n)}$,
we arrive at the following

\begin{definition}
  The free quasi-symmetric generating function $\bGamma(P)\in\FQSym$
  of a poset $P$ is
  \begin{equation}
    \bGamma(P)=\sum_{\sigma\in L(P)} \F_\sigma\,.
  \end{equation}
\end{definition}
This amounts to encode a poset by the set of its linear extensions. It
is well known, and easy to see, that if $P_1$ is an order on $[k]$ and
$P_2$ an order on $[l]$, the order $P=P_1\sqcup P_2$ on $[n]=[k+l]$
defined by $i<_P j \Leftrightarrow i<_{P_1} j$ or $i-k <_{P_2} j-k$
has for linear extensions the shifted shuffles of those of $P_1$ and
$P_2$:
\begin{equation}
  \sum_{\sigma\in L(P)}\sigma=\sum_{\alpha\in L(P_1)}\sum_{\beta\in L(P_2)}
  \alpha\shuffle \beta[k]\,.
\end{equation}
Thus, for the free generating functions, one has as well
\begin{equation}
  \bGamma(P_1\sqcup P_2) = \bGamma(P_1)\bGamma(P_2)
\end{equation}
in $\Z\free{A}$.

It will be convenient to introduce the notation $P_2[k]$ for the order
on $[k+1,k+l]$ defined above, so that $L(P_1\sqcup P_2)=L(P_1)\shuffle
L(P_2[k])$.

\begin{example}
  The free Schur functions $\SS_t$ are of the form $\bGamma(P)$ for
  the posets associated to plane partitions. Malvenuto \cite{Mv} has shown that
  if $\bGamma(P)\in\FSym$, then $P$ is associated to a plane
  partition. This is a step towards a famous conjecture of Stanley,
  asserting that the conclusion remains valid as soon as the
  commutative image of $\bGamma(P)$ is symmetric.
\end{example}

\begin{example}
  The concatenation $P_1\sqcup P_2$ is not the only interesting poset
  which can be constructed from $P_1$ and $P_2$.  One can also define
  $P=P_1\wedge P_2$ as the poset obtained by adjoining a maximal
  element to the juxtaposition of $P_1$ and $P_2$.  The correct way to
  do this is to take as maximal element $h=k+1$ if $P_1$ is a poset on
  $[k]$. Therefore, $i<_P j$ iff $i,j \le k$ and $i<_{P_1} j$ or
  $i,j>k+1$ and $i-k-1<_{P_2}i-k-1$, or $j=k+1$. The linear extensions
  of $P$ are clearly
  \begin{equation}
    L(P_1\wedge P_2)= (L(P_1)\shuffle L(P_2[h]))\cdot h\,.
  \end{equation}

The posets generated from $\bullet=[1]$ by the operation $\wedge$ are in
one-to-one correspondence with binary trees, since they correspond
to all possible bracketings of the words $\bullet\bullet\cdots\bullet$.
Let $\F(T)=\bGamma(T)$ be the free quasi-symmetric generating functions
of such posets. We will see that they span a subalgebra of $\FQSym$,
which is precisely the Hopf algebra of binary trees introduced by
Loday and Ronco \cite{LR}. Indeed, let $T=T_1\wedge T_2$
and $T'=T'_1\wedge T'_2$ be two binary trees. From the above
considerations,
we see that
\begin{multline*}
  L(T\sqcup T') =L(T)\shuffle L(T'[n])
  = 
  \sum_{\substack{\alpha\in L(T_1) \\ \beta\in L(T_2)}}\ 
  \sum_{\substack{\alpha'\in L(T'_1[n]) \\ \beta'\in L(T'_2[n])}}
  \left[ (\alpha\shuffle \beta[h])h\right]
  \shuffle
  \left[
    (\alpha'\shuffle\beta'[h'])h'
  \right] \\
  = \sum
  \left[
    (\alpha\shuffle\beta[h]) \shuffle (\alpha'\shuffle\beta'[h'])h'
  \right] h
  + \sum
  \left[
    ((\alpha\shuffle\beta[h])h) \shuffle (\alpha'\shuffle \beta'[h'])
  \right] h'
\end{multline*}
(using the formula $(ua)\shuffle (vb) = (u\shuffle vb)a+ (ua\shuffle v)b$,
valid for $a,b\in A$ and $u,v\in A^*$). Therefore,
\begin{equation}
  L(T\sqcup T')
  =  L(T_1\wedge(T_2\sqcup T')) +L ((T\sqcup T'_1)\wedge T_2)
\end{equation}
which proves that $\F(T)\F(T')$ is a sum of elements $\F(T'')$ which
are given by the above recursion.

The connection with the algebra of Loday and Ronco comes for the fact
that $\F(T)=\sum_{\T(\sigma)=T}\G_\sigma$, where $\T(\sigma)$ is the
underlying binary tree of $\sigma$, defined as follows: if $n=1$
($\sigma$ is the empty word), $\T(\sigma)=\bullet$, otherwise, write
$\sigma=unv$, $\alpha=\std(u)$, $\beta=\std(v)$. Then,
$\T(\sigma)=\T(\alpha)\wedge \T(\beta)$, where $T_1\wedge T_2$ is the
binary tree having $T_1$ as left subtree and $T_2$ as right subtree.
\end{example}

\subsection{Posets as $0$-Hecke modules}

There is a striking similarity between the behavior of the
quasi-symmetric generating functions of posets under concatenation,
and the characteristic quasi-symmetric functions of $0$-Hecke modules
under induction product. Actually, the former is a special case of the
latter:

\begin{definition}
  The \emph{$0$-Hecke module $M_P$ associated with a poset $P$} is the
  (right) $0$-Hecke module with basis the set of linear extensions
  $L(P)$ and structure defined by
  \begin{equation}
    \sigma T_i =
    \left\{
      \begin{array}{cl}
        \sigma\sigma_i & \text{if $i\not\in \Des(\sigma)$ and
                               $\sigma\sigma_i\in L(P)$\,,} \\[1mm]
        0              & \text{if $i\not\in \Des(\sigma)$ but
                                $\sigma\sigma_i\not\in L(P)$\,,} \\[1mm]
        -\sigma        & \text{if $i\in \Des(\sigma)$\,.}
    \end{array}
    \right.
  \end{equation}
\end{definition}

\Proof We have to prove that $M_P$ is actually a $0$-Hecke module.
Here we need some definitions.
\begin{definition}
  A poset $P$ is \emph{rise free} if there is no $i<j$ such that
  $i<_Pj$. 
\end{definition}
Recall that each poset has a minimal linear extension $E(P)$ defined by
\begin{equation}
  i<_{E(P)} j
  \qquad\text{iff}\qquad
  i<_P j\text{ or } (i<j\text{ and } j\not<_P i)\,.
\end{equation}
One easily has 
\begin{proposition}
  Let $P$ a rise free poset. The set of permutations that are larger
  for the right weak order than $E_P$ is exactly the set of linear
  extensions of $P$.
\end{proposition}

It has for consequence that if $P$ is a rise free poset, the submodule
of the regular representation generated by its minimal linear
extension has the structure defined above. Then $M_P$ is a module
for $P$ rise free.

\newcommand{\RF}{\operatorname{RiseFree}}
Now, if $P$ is not rise
free, let $\RF(P)$ be its associated rise free poset defined by
\begin{equation}
  i<_{\RF(P)}j\qquad\text{iff}\qquad (i>j\text{ and } i<_P j)\,.
\end{equation}
Note that $P$ and $\RF(P)$ have the same minimal linear extension.
Consider the module $M_{\RF(P)}$. It has for basis the set of
permutations that are greater than the minimal linear extension of
$\RF(P)$ and $P$. If a permutation $\sigma$ in this set is not a
linear extension of $P$ then there is a $i<j$ such that $i<_P j$ but
$\sigma_i>\sigma_j$. And then all the permutations bigger than
$\sigma$ are not linear extension of $P$. This means that the set of
permutations larger than $E(P)$ but that are not linear extensions of
$P$ span a sub-module $N$ of $M_{\RF(P)}$. Now it is easy to see that
\begin{equation}
  M_P \equiv M_{\RF(P)}/ N\,.
\end{equation}
is a realisation of $M_P$. And hence $M_P$ is a module.\qed
\bigskip

\begin{figure}[ht]
  \vskip1cm
  \epsfig{file=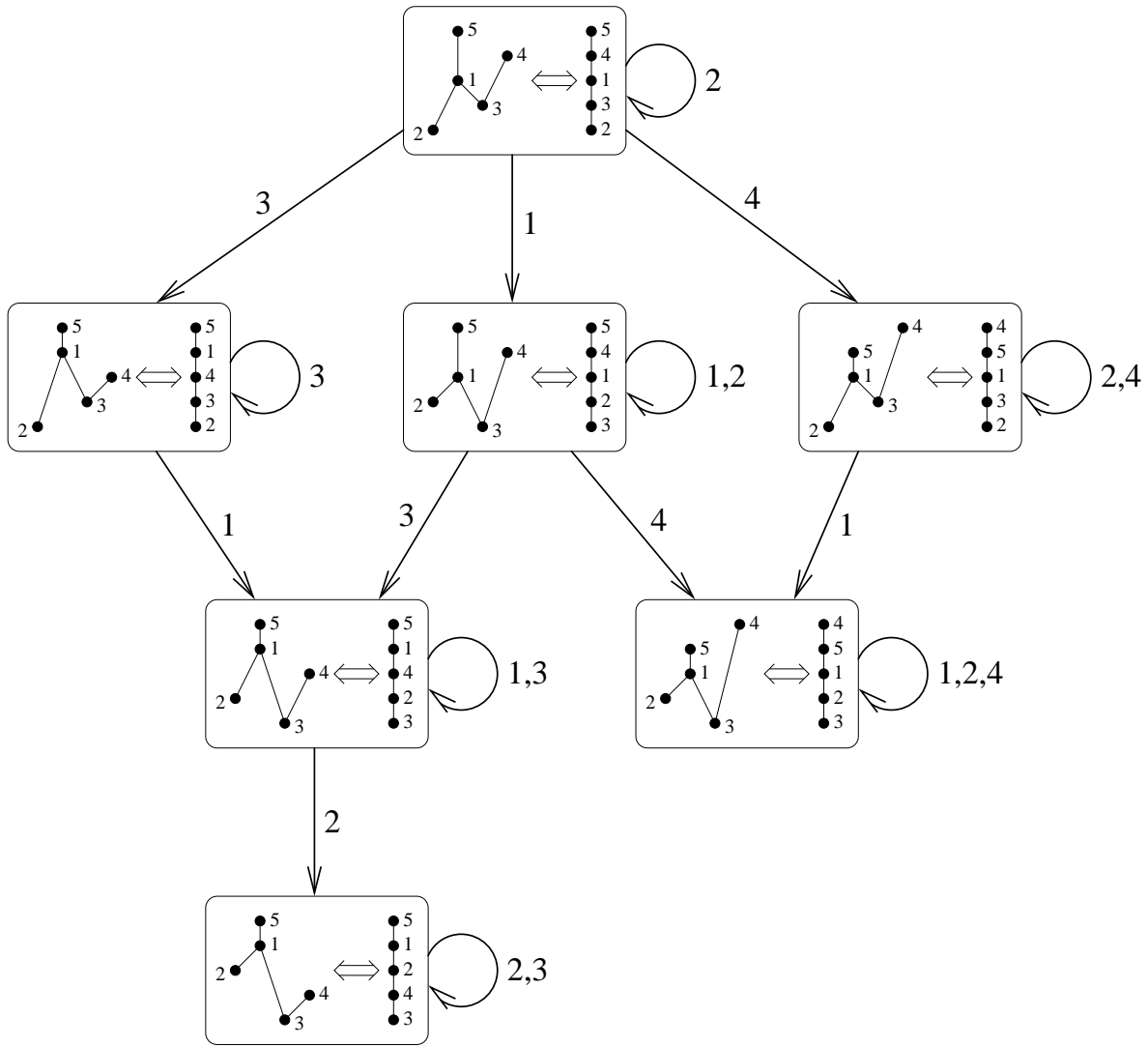,width=\textwidth}
  \caption{\label{fig::Posmod}Example of module associated with a
    poset}
  \vskip1cm
  To each vertex of the graph enclosed in a box corresponds a basis
  element associated with the depicted linear extension. There is a
  straigth arrow labelled $i$ from $u$ to $v$ if $uT_i=v$. A loop
  labelled $i$ means that $uT_i=-u$.  If there is no arrow
  labelled $i$ leaving the vertex $u$, then $uT_i=0$.
\end{figure}

Then one has following proposition. 
\begin{proposition}
  \begin{itemize}
  \item[(i)] Let $P$ be a poset. Then  $\ch(M_P)=\Gamma(P)$. 
  \item[(ii)] $M_{P\sqcup P'} = M_P \hat\otimes M_{P'}$ and
  consequently 
  \begin{equation}
     \fch(M_{P\sqcup P'}) = \bGamma(P) \,\bGamma(P')
                          = \fch(M_{P})\,\fch(M_{P'})
                          = \fch(M_P \hat\otimes M_{P'})
  \end{equation}
  \end{itemize}
\end{proposition}
See Figure \ref{fig::Posmod} for an example.

\newpage
\renewcommand{\topfraction}{1}  
\renewcommand{\textfraction}{0} %

\newpage

\subsection{Shuffle and pseudo-convolution}

In this section, we will encounter another example of
pseudo-convolution, in the sense of \ref{sec:conv}. Let $\Box$ be the
operation on $\bigoplus \C\SG_n$ defined on permutations by
\begin{equation}
  \alpha\Box \beta =
  \sum_{I\sqcup J = [1,k+l]}
     i_{\alpha(1)}\cdots i_{\alpha(k)}
     \shuffle
     j_{\beta(1)}\cdots j_{\beta(l)}
\end{equation}
This operation arises naturally in the problem of calculating the
orthogonal projection onto the free Lie algebra.  One can show that
this problem boils down to the inversion of the element
\begin{equation}
  T_n = \sum_{k=0}^n (1\cdots k)\Box (1\cdots n-k)
\end{equation}
of $\Q\SG_n$ \cite{Du,Kly}. No closed formula is known for $T_n^{-1}$,
but numerical experiments suggest that it should be possible to give a
combinatorial description of its characteristic polynomial.

\begin{example} The characteristic polynomial of $T_4$ as an operator on
  the regular representation of $\SG_4$ is
  $$
  (x-2)^6(x-6)^4(x-14)^3(x-18)^3(x-42)^3(x-70)(x^2-28x+84)^2\,.
  $$
\end{example}

It is natural to introduce $q$-{\it pseudo-convolution} $\Box_q$,
which is defined similarly, with $\shuffle$ replaced by $\shuffle_q$.
Actually, this operation can be interpreted in $\FQSym$ in the same
way as the ordinary convolution. Let
$\pair{\sigma}{\tau}=\delta_{\sigma,\tau}$ be the scalar product on
the group algebra for which permutations form an orthonormal basis, so
that
\begin{equation}
  \alpha \Box_q \beta = \sum_\sigma \pair{\sigma}{\alpha\Box_q\beta}\sigma\,.
\end{equation}

\begin{proposition}
  The algebra of free quasi-symmetric functions is a $q$-shuffle
  subalgebra of $\C(q)\free{A}$, and in the $\G$-basis, the structure
  constants coincide with those of $q$-pseudo-convolution
  \begin{equation}
    \G_\alpha \shuffle_q \G_\beta
    = \sum_\sigma \pair{\sigma}{\alpha\Box_q\beta}\G_\sigma \,.
  \end{equation}
\end{proposition}

\Proof (sketch) We proceed as for Proposition \ref{MULTF}. In the
biword notation, we have
\begin{align*}
  \G_\alpha \shuffle_q \G_\beta
  &= \sum_{\substack{u\uparrow\alpha^{-1}\\ v\uparrow\beta^{-1}}}
         \binom{u}{\alpha^{-1}} \shuffle_q \binom{v}{\beta^{-1}} 
    = \sum_{\substack{u\uparrow\alpha^{-1}\\ v\uparrow\beta^{-1}}} 
         \binom{uv}{(\alpha\shuffle_q \beta)^{\vee}}
  \\
  &= \sum_{\sigma\in\alpha\Box_q \beta}
         \sum_{w\uparrow\sigma^{-1}}
             \binom{w}{\sigma}
    = \sum_\sigma \pair{\sigma}{\alpha\Box_q\beta} \G_\sigma\,.
\end{align*}
  
In particular, for $q=1$, identifying $\G_{(12\cdots n)}$ to the
noncommutative complete function $S_n$, we see that
\begin{equation}
  T_n = S_n + S_1\shuffle S_{n-1} + \cdots +S_{n-1}\shuffle S_1 +S_n
      = h_n(2X)
\end{equation}
if we identify the $\shuffle$-subalgebra generated by the $S_n$
with the algebra of commutative symmetric functions of some alphabet
$X$. At this point, it is natural to introduce $q$-analogues $T_n(q)$.
If we define them as
\begin{equation}
  T_n(q) = S_n +q S_1\shuffle_q S_{n-1}+\cdots+q^n S_n
\end{equation}
we see that $T_n(q)=S_n((1+q)B)$, if we now identify the
$\shuffle_q$-subalgebra of $\FQSym$ generated by the $S_n$ with
$\Sym(B_q)$ for a noncommutative alphabet $B_q$. That is, we have a
one-parameter family of identifications of the $S^I$ with elements of
the group algebra.

\begin{example}
The characteristic polynomial of $T_3(q)$ is
$$
(x-2)^2(x-4-4q-2q^2)^2(x-8-6q-6q^2)(x-4+2q-2q^2)\,.
$$
\end{example}

It makes sense to consider the quasi-symmetric generating functions of
the elements $T_n$, which amounts to take the commutative images of
the corresponding elements of $\FQSym$ (here it does not matter
whether one interprets $\sigma$ as $\F_\sigma$ or $\G_\sigma$ since
$T_n$ is self-adjoint. One finds that
\begin{equation}
  \underline T_n = \sum_{i+j=n}\binom{n}{i}h_ih_j
  = [t^n] \frac{1}{1-t}\sum_{m=0}^n \left( \frac{t}{1-t}\right)^m h_mh_{n-m}\,.
\end{equation}
The first values are, on the Schur basis
$$
2s_1,\ 4s_2+2s_{11},\ 8s_3+6s_{21},\ 16 s_4 +14 s_{31} + 6s_{22}\, \dots
$$
The elementary symmetric functions of the $\shuffle$-algebra
generated by the $S_n=\G_{(12\cdots n)}$ also seem to be interesting.
It would be interesting to investigate the structure of $\FQSym$ as a
$\shuffle$-module over this commutative subalgebra, and also the
$q$-analogue of this situation.

This suggests the possibility of using the machinery of noncommutative
symmetric functions to invert $T_n(q)$. The problem is to interpret
the internal product of $\Sym(B_q)$ in terms of the structure of
$\FQSym$, and more precisely to connect it to the ordinary composition
of permutations. That is, if one defines $*_q$ on $\Sym(B_q)$ by the
standard formulas giving $S^I*S^J$, for example, does there exist an
automorphism $\phi_q$ of $\C(q)\free{A}$ such that
$F*G=\phi_q^{-1}(\phi_q(F)\circ \phi_q(G))$? (here $\circ$ is the
composition of permutations).

\subsection{Identities}

A few identities between series of free quasi-symmetric functions
(mainly conjectures) can be found in \cite{U}. For example, the
inverses of the series
\begin{eqnarray*}
  H_1 &=& \sum_I (-1)^{\ell(I)}\F_{\omega(I)}\\
  H_2 &=& \sum_{n\ge 0} (-1)^n\F_{\omega(2^n)}\\
  H_3 &=&  \sum_I (-1)^{\ell(I)}\F_{\omega(2I)}\\
\end{eqnarray*}
are conjectured to be as follows. For a permutation $\sigma$ of shape
$I$, let $\hat\sigma=\sigma\alpha(I)$. Then,
\begin{eqnarray*}
  H_1^{-1} &=& \sum_\alpha \G_{\hat\alpha}\\
  H_2^{-1} &=& \sum_\beta \G_{\hat\beta}\\
  H_3^{-1}&=&\sum_\gamma \G_{\hat\gamma}
\end{eqnarray*}
where $\alpha$ runs over all permutations, $\beta\in\SG_{2p}$ runs
over permutations of shape $2^{2p}$, and $\gamma\in\SG_{2p}$ runs over
permutations with descent set contained in $\{2,4,\ldots,2p-2\}$.

                   \section{The $0$-Hecke algebra revisited}

\subsection{$H_n(0)$ as a Frobenius algebra}

Recall that a bilinear form $(\,,\,)$ on a $\K$-algebra $A$ is said to
be associative if $(ab,c)=(a,bc)$ for all $a,b,c\in A$, and that $A$
is called a {\em Frobenius algebra} whenever it has a nondegenerate
associative bilinear form. Such a form induces an isomorphism of left
$A$-modules between $A$ and the dual $A^*$ of the right regular
representation. Frobenius algebras are in particular self-injective,
so that finitely generated projective and injective modules coincide
(see \cite{CR}).

For a basis $(Y_\sigma)$ of $H_n(0)$, we denote by $(Y^*_\sigma)$ the
dual basis. We set $\chi=T^*_\omega$, where $\omega=(n\, n-1\ldots 1)$
is the longest permutation of $\SG_n$.

\begin{proposition}
  \begin{itemize}
  \item[(i)] The associative bilinear form defined by
    \begin{equation}
      (f,g)=\chi(fg)
    \end{equation}
is non-degenerate on $H_n(0)$. Therefore, $H_n(0)$ is a
    Frobenius algebra.
  \item[(ii)] $(\eta_\sigma,\eta_{\tau^{-1}\omega})=\delta(\sigma\ge
    \tau)$, where $\ge$ is the Bruhat order on $\SG_n$, and for a
    statement $P$, $\delta(P)$ is 1 when $P$ is true and 0 otherwise.
  \item[(iii)] The elements
    $\zeta_\sigma=(-1)^{\ell(\omega\sigma^{-1})}\xi_{\omega\sigma^{-1}}$
    satisfy
    \begin{equation}
      (\zeta_\sigma,\eta_\tau)=\delta_{\sigma,\tau}\,.
    \end{equation}
  \end{itemize}
\end{proposition}

\Proof
The bilinear form defined in {\it (i)} is clearly associative. That it is
non-degenerate follows from { \it (ii)}, which implies that the matrix
$(\eta_\sigma,\eta_\tau)$ is, up to a permutation of columns, the incidence
matrix of the Bruhat order, which is obviously invertible. 
The proof of {\it (ii)} is a simple induction on $\ell(\sigma)$.
Finally, {\it (iii)} follows from {\it (ii)} and \cite{La}, Lemme 1.13,
which says that
\begin{equation}
\xi_\alpha =\sum_{\beta\le\alpha}T_\beta\quad
{\rm and}
\quad
\eta_{\alpha}=\sum_{\beta\le\alpha}(-1)^{\ell(\beta)}\xi_\beta\,.
\end{equation}

\bigskip

\begin{remark}
  As recently shown by L. Abrams \cite{Ab}, a Frobenius algebra is endowed
  with a comultiplication $\delta:\ A\rightarrow A\otimes A$
  which is a morphism of $A$-bimodules, that is, $\delta(axb)=a\delta(x)b$.
  It can be defined by the formula
  \begin{equation}
    \delta=(\lambda^{-1}\otimes \lambda^{-1})\circ (\mu\circ T)^*\circ \lambda
  \end{equation}
  where $\lambda:\ A\rightarrow A^*$ is an isomorphism of left $A$-modules,
  $\mu:\ a\otimes b \mapsto ab$ is the multiplication map,
  and $T:\ a\otimes b\mapsto b\otimes a$ is the exchange operator.
  Since $\delta$ is a bimodule map, it is completely specified
  by the element $\delta(1_A)$, which we will now calculate
  explicitly for $H_n(0)$.
  Let $\lambda$ be defined by $\lambda(x)(y)=(y,x)$. Then,
  $$
  \lambda(\eta_\sigma)=\sum_{\omega\tau^{-1}\le \sigma}\eta^*_\tau
  $$
  so that $\lambda^{-1}(\eta^*_\sigma)=\zeta_\sigma$.
  If we define the permutation $\{\alpha,\beta\}$ by the rule
  $\eta_\alpha\eta_\beta=\eta_{\{\alpha,\beta\}}$, then,
  \begin{align*}
    \delta(1)
    & = \sum_{\{\alpha,\beta\}=\omega}\zeta_\beta\otimes\zeta_\alpha \\
    & = \sum_\alpha
          \left(
            \sum_{\{\alpha,\beta\}=\omega}
               (-1)^{\ell(\omega\beta^{-1})} \xi_{\omega\beta^{-1}}
          \right) \otimes \zeta_\alpha \\
    & = \sum_\alpha
           \left(
              \sum_{\gamma\le\alpha}
                 (-1)^{\ell(\gamma)}\xi_\gamma
           \right) \otimes \zeta_\alpha
      = \sum_\alpha \eta_\alpha\otimes \zeta_\alpha\,.
  \end{align*}
  Therefore, the canonical comultiplication of $H_n(0)$ is given by
  \begin{equation}
    \delta(1)=\sum_{\sigma\in\SG_n} \eta_\sigma\otimes\zeta_\sigma
  \end{equation}
  and $\delta(x)=x\delta(1)=\delta(1)x$.
\end{remark}

\subsection{$\FQSym$ as a Grothendieck ring}

Let $(g_\sigma)$ be the basis of $H_n(0)$ defined by
\begin{equation}
  g_\sigma = T_{\sigma\alpha(I)^{-1}} \epsilon_I
\end{equation}
where $I=C(\sigma)$ is the descent composition of $\sigma$ and
$\epsilon_I$ the  generator of the principal indecomposable
projective module $\P_I$. Then, $\{g_\sigma|\sigma\in[\alpha(I),\omega(I)]\}$
is a basis of $\P_I$ (the interval is taken with respect to the weak order). 

\begin{definition}
  For any permutation $\sigma\in\SG_n$, we denote by $\NN_\sigma$ the
  submodule of $\P_I$ (where $I=C(\sigma)$) generated by $g_\sigma$.
\end{definition}

All the $\NN_\sigma$ are indecomposable $H_n(0)$-modules, since any
submodule of a $\P_I$ must contain its one-dimensional socle, and
therefore cannot be a direct summand. The simple $H_n(0)$ modules are
the $\NN_{\omega(I)}$, and $\P_I=\NN_{\alpha(I)}$.

Of course, the $\NN_\sigma$ do not exhaust all submodules of the
$\P_I$, but, as we will see, they generate an interesting subcategory
${\cal N}_n$ of $H_n(0)-\mod$. In particular, all the specializations
$q=0$ of the Specht modules $V_\lambda(q)$ of $H_n(q)$, as well as
their  skew versions $V_{\lambda/\mu}(q)$, with
$\lambda/\mu$ connected, are of the form $\NN_\sigma$, where $\sigma$
is the row reading of the hyperstandard tableau of shape $\lambda$ (or
$\lambda/\mu$) i.e., the tableau whose columns are filled with
consecutive integers. As a consequence, all the $V_{\lambda/\mu}(0)$,
with $\lambda/\mu$ connected, are indecomposable.

Define a characteristic map with values in $\FQSym$ by
\begin{equation}
  \fch(\NN_\sigma) = N_\sigma=\sum_{\tau\in[\sigma,\omega(I)]}\G_\sigma\,.
\end{equation}
This definition is compatible with the former one for projective
modules, since $\fch(\P_I)=R_I$.  More generally, the characteristic
of a Specht module is a free symmetric function:
$\fch(V_\lambda(0))=\SS_t$, where $t$ is the tableau congruent to the
contretableau of shape $\omega(\lambda)$ whose rows consist of
consecutive integers (e.g., 456 23 1 for $\lambda=(321)$).

\begin{proposition}
  The characteristic map is compatible with induction
  product, that is, we have an exact sequence
  \begin{equation}
    0 \rightarrow \NN_\beta
      \rightarrow \NN_\sigma\hat\otimes\NN_\tau \rightarrow \NN_\alpha 
      \rightarrow 0
  \end{equation}
  where $\alpha=\sigma\bullet\tau$, and if as words
  $\sigma^{-1}=ukv$, $\tau^{-1}[k]=u'(k+1)v'$ then
  $\beta^{-1}=uu'(k+1)kvv'$, and also
  \begin{equation}
    \fch(\NN_\sigma\hat\otimes\NN_\tau )
       = N_\sigma N_\tau = N_\alpha + N_\beta\,.
  \end{equation}
  In the case of skew Specht modules indexed by connected skew diagrams
  $D,D'$, the formula reads
  \begin{equation}
    0 \rightarrow V_{D_2} \rightarrow V_D\hat\otimes V_{D'}
      \rightarrow V_{D_1} \rightarrow 0
  \end{equation}
  where $D_1$ and $D_2$ are the two ways of glueing the first box
  of $D'$ to the last box of $D$.
\end{proposition}

\Proof Remark first that if $C(\sigma)=I$ and $C(\tau)=J$,
$M=\NN_\sigma\hat\otimes\NN_\tau $ is a submodule of
$\P_I\hat\otimes\P_J=\P_{I\triangleright J}\oplus \P_{I\cdot J}$.
Also, $M$ is a combinatorial module. It is generated by the element
$g_\sigma\otimes g_\tau$, which can be represented by the skew ribbon
$r_0$ obtained by making the upper left corner of the first cell of
the ribbon of $\tau[k]$ coincide with the bottom right corner of the
last cell of the ribbon of $\sigma$. The combinatorial basis of $M$ is
formed by those skew ribbon of the same shape as $r_0$ which can be
obtained from $r_0$ by application of a chain of operators
$\eta_i=-T_i$. Their action is given by the same formulas as for the
case of connected ribbons representing the bases of the projective
indecomposable modules: if $i$ is a recoil of $r$, then $\eta_i(r)=r$.
If $i+1$ is in the same row as $i$, then $\eta_i(r)=0$, and otherwise,
$\eta_i(r)=r'$, the skew ribbon obtained from $r$ by exchanging $i$
and $i+1$.

Now, the skew ribbons generated from $r_0$ can be converted into
connected ribbons of shape $IJ$ or $I\triangleright J$, according to
whether the first entry of the right connected component is greater or
smaller than the last entry of the left component. The generator $r_0$
corresponds to the shape $I\triangleright J$, filled with the
permutation $\alpha$. According to the above rules, the action of
$H_n(0)$ will generate all permutations of this shape which are
greater than $\alpha$ for the weak order, plus some other ones of
shape $IJ$.

All the permutations of shape $IJ$ are greater than those of shape
$I\triangleright J$, and span therefore a submodule, which is easily
seen to be generated by $\beta$.  Indeed, define $\beta$ as the
smallest (for the weak order) permutation of shape $IJ$ which is
greater than $\alpha$.  Set $\beta=st$ as a word, with $|s|=k$. Since
$\beta >\alpha$, we have $\std(s)=\sigma$ and $\std(t)=\tau$. This
means that the letters $1,\ldots,k$ occur in the same order in
$\sigma^{-1}$ and in $\beta^{-1}$, and also, $k+1,\ldots,k+l$ occur in
the same order in $\tau^{-1}$ and $\beta^{-1}$. Hence, $\beta^{-1}\in
\sigma^{-1}\shuffle \tau^{-1}[k]$. Also, $k$ must be a descent of
$\beta$. Hence, in $\beta^{-1}$, the letter $k+1$ appears on the left
of $k$. The smallest permutation with these properties is
$\beta^{-1}=uu'(k+1)kvv'$, as claimed.

Hence, $\NN_\beta$ is a submodule (even a subgraph) of $M$, and the
quotient is isomorphic to $\NN_\alpha$. Now, the permutations obtained
by applying the $\eta_i$ to $r_0$ can also be described as those
$\gamma$ which, as words, satisfy $\gamma=uv$ with
$\std(u)\in[\sigma,\omega(I)]$ and $\std(v)\in[\tau,\omega(J)]$. These
are exactly the standardizations of the words occurring in the product
$N_\sigma N_\tau$.
\qed
\bigskip\pagebreak[3]

In particular, we obtain a description of the induction products
of simple modules, which is much more precise than the one given
by the product of quasi-symmetric functions:
\begin{corollary}
  Any induction product of simple modules
  $\SS_{I_1}\hat\otimes\cdots\hat\otimes \SS_{I_r}$ has a filtration
  by modules $\NN_\sigma$, which can be explicitely computed.
\end{corollary}

By using a standard result on self-injective algebras, we can now
define another family of indecomposable modules. Indeed, for any
self-injective Artin algebra $A$, and any exact sequence
\begin{equation}
  0\rightarrow N\rightarrow P\rightarrow M\rightarrow 0\,,
\end{equation}
of left $A$-modules, with $P$ projective, $N$ is indecomposable non
injective, and $N\rightarrow P$ an injective hull, iff $M$ is
indecomposable non projective, and $P\rightarrow M$ a projective
cover ({\it cf. } \cite{CR}).  It is customary to set $N=\Omega M$ and $M=\Omega^{-1}N$.
$\Omega$ is called the syzygy functor (as defined here it is only a
map on the set of modules, but it becomes a functor in the stable
category; here it is well defined as a map because of the unicity of
the minimal projective resolution).

Since the inclusion $\NN_\sigma\rightarrow \P_I$ is clearly an
injective hull, we have:

\begin{lemma}
  For $\sigma\in]\alpha(I),\omega(I)]$, $\MM_\sigma=\P_I/\NN_\sigma$
  is indecomposable.
\end{lemma}

Starting with $M$ simple, next taking a projective cover of $N$, and
iterating the process, one can construct a sequence of indecomposable
modules $\Omega^nM$. In this way, one can see that for $n>3$, $H_n(0)$
is not representation finite: the sequences $\Omega^n\SS_I$ are
neither finite nor periodic for $I\not=(n),(1^n)$, and $\dim_\C
\Omega^n\SS_I \rightarrow \infty$.

\subsection{Homological properties of $H_n(0)$ for small $n$}

Being a finite dimensional elementary $\C$-algebra, $H_n(0)$ can be
presented in the form $\C Q/{\cal I}$, where $Q$ is a quiver, $\C Q$ its
path algebra, and $\cal I$ an ideal contained in ${\cal J}^2$ where
$\cal J$ is the ideal generated by all the arrows of $Q$
\cite{ARS}. The vertices of $Q$ are the simple modules $\SS_I$, and
the number $e_{IJ}$ of arrows $\SS_I\rightarrow \SS_J$ is equal to
$\dim \Ext^1(\SS_I,\SS_J)=[\rad\P_I/\rad^2\P_I : \SS_J]$.

Therefore, $e_{IJ}=c_{IJ}^{(1)}$, where
$c_{IJ}^{(k)}=[\rad^k\P_I/\rad^{k+1}\P_I : \SS_J]$ are the
coefficients of the $q$-Cartan invariants
\begin{equation}
  c_{IJ}(q) =\sum_{k\ge 0} c_{IJ}^{(k)} q^k
\end{equation}
associated to the radical series. Let
$C_n(q)=(c_{IJ}(q))_{I,J\models n}$. For $n\le 4$, these matrices are
as follows.

$$
\begin{array}{|c||c|c|c|c|}
\hline
& 3 & 21 & 12 & 1111\\
\hline\hline
3 & 1 & 0 & 0 & 0 \\
\hline
21& 0 & 1 & q & 0 \\
\hline
12& 0 & q & 1 & 0\\
\hline
111&0&0&0&1\\
\hline
\end{array}
$$

\bigskip 

$$
\begin{array}{|c||c|c|c|c|c|c|c|c|}
\hline
& 4 & 31 & 22 & 211 & 13 & 121 & 112 & 1111\\
\hline\hline
4 & 1&0&0&0&0&0&0&0\\
\hline
31 &0& 1&q&0&q^2&0&0&0\\
\hline
22&0&q&1+q^2&0&q&q&0&0\\
\hline
211 &0 & 0&0&1&0&q&q^2&0\\
\hline
13 & 0 &q^2&q&0&1&0&0&0\\
\hline
121&0&0&q&q&0&1+q^2&q&0\\
\hline
112&0&0&0&q^2&0&q&1&0\\
\hline
1111&0&0&0&0&0&0&0&1\\
\hline
\end{array}
 $$                             

\begin{sidewaystable}
\centering
\begin {tabular}{|c||c|c|c|c|c|c|c|c|c|c|c|c|c|c|c|c|} 
\hline
  & $ 5 $ & $ 41 $ & $ 32 $ & $ 311 $ & $ 23 $ & $ 221 $ & $ 212 $ & $ 2111 $ & $ 14 $ & $ 131 $ & $ 122 $ & $ 1211 $ & $ 113 $ & $ 1121 $ & $ 1112 $ & $ 11111$ \\
\hline\hline
5  & $ 1$ & $0$ & $0$ & $0$ & $0$ & $0$ & $0$ & $0$ & $0$ & $0$ & $0$ & $0$ & $0$ & $0$ & $0$ & $0$\\
\hline
41  & $0$ & $  1$ & $q$ & $0$ & ${q}^{2}$ & $0$ & $0$ & $0$ & ${q}^{3}$ & $0$ & $0$ & $0 $ & $0$ & $0$ & $0$ & $0$\\
\hline
32  & $0$ & $ q$ & $1+{q}^{2}$ & $0$ & ${q}^{4}+q$ & $q$ & $0$ & $0$ & ${q}^{2}$ & ${q}^{2}$ & ${q}^{3}$ & $0$ & $0$ & $0$ & $0$ & $0$\\
\hline
311  & $0$ & $ 0$ & $0$ & $1$ & $0$ & $q$ & ${q}^{2}$ & $0$ & $0$ & ${q}^{2}$ & ${q}^{3}$ & $0$ & ${q}^{4}$ & $0$ & $0$ & $0$\\
\hline
23  & $0$ & $ {q}^{2}$ & ${q}^{4}+q$ & $0$ & $1+{q}^{2}$ & ${q}^{3}$ & $0$ & $0$ & $q$ & ${q}^{2}$ & $q$ & $0$ & $0$ & $0$ & $0$ & $0$\\
\hline
221  & $0$ & $ 0$ & $q$ & $q$ & ${q}^{3}$ & $2\,{q}^{2}+1$ & $q+{q}^{3}$ & $0$ & $0$ & $q+{q}^{3}$ & ${q}^{4}+2\,{q}^{2}$ & $q$ & ${q}^{3}$ & ${q}^{3}$ & $0 $ & $0$\\
\hline
212  & $0$ & $0$ & $0$ & ${q}^{2}$ & $0$ & $q+{q}^{3}$ & $1+{q}^{4}$ & $0$ & $0$ & ${q}^{2}$ & $q+{q}^{3}$ & ${q}^{2}$ & ${q}^{2}$ & ${q}^{2}$ & $0$ & $0$\\
\hline
2111  & $0$ & $ 0$ & $0$ & $0$ & $0$ & $0$ & $0$ & $1$ & $0$ & $0$ & $0 $ & $q$ & $0$ & ${q}^{2}$ & ${q}^{3}$ & $0$\\
\hline
14  & $0$ & ${q}^{3}$ & ${q}^{2}$ & $0$ & $q$ & $0$ & $0$ & $0$ & $1$ & $0$ & $0 $ & $0$ & $0$ & $0$ & $0$ & $0$\\
\hline
131  & $0$ & $ 0$ & ${q}^{2}$ & ${q}^{2}$ & ${q}^{2}$ & $q+{q}^{3}$ & ${q}^{2}$ & $0$ & $0$ & $1+{q}^{4}$ & $q+{q}^{3}$ & $0$ & ${q}^{2}$ & $0$ & $0$ & $0$\\
\hline
122  & $0$ & $ 0$ & ${q}^{3}$ & ${q}^{3}$ & $q$ & ${q}^{4}+2\,{q}^{2}$ & $q+{q}^{3}$ & $0$ & $0$ & $q+{q}^{3}$ & $2\,{q}^{2}+1$ & ${q}^{3}$ & $q$ & $q$ & $0$ & $0$\\
\hline
1211  & $0$ & $ 0$ & $0$ & $0$ & $0$ & $q$ & ${q}^{2}$ & $q$ & $0$ & $0$ & ${q}^{3}$ & $1+{q}^{2}$ & $0 $ & ${q}^{4}+q$ & ${q}^{2}$ & $0$\\
\hline
113  & $0$ & $  0$ & $0$ & ${q}^{4}$ & $0$ & ${q}^{3}$ & ${q}^{2}$ & $0$ & $0 $ & ${q}^{2}$ & $q$ & $0$ & $1$ & $0$ & $0$ & $0$\\
\hline
1121  & $0$ & $ 0$ & $0$ & $0$ & $0$ & ${q}^{3}$ & ${q}^{2}$ & ${q}^{2}$ & $0 $ & $0$ & $q$ & ${q}^{4}+q$ & $0$ & $1+{q}^{2}$ & $q$ & $0$\\
\hline
1112  & $ 0$ & $0$ & $0$ & $0$ & $0$ & $0$ & $0$ & ${q}^{3}$ & $0$ & $0 $ & $0$ & ${q}^{2}$ & $0$ & $q$ & $1$ & $0$\\
\hline
11111  & $ 0$ & $0$ & $0$ & $0$ & $0$ & $0$ & $0$ & $0$ & $0$ & $0$ & $0$ & $0$ & $0$ & $0$ & $0$ & $1$\\
\hline
\end {tabular}
\caption{The $q$-Cartan matrix of $H_5(0)$}
\end{sidewaystable}


\noindent The corresponding quivers are given on Figures
\ref{fig:quiv34} and \ref{fig:quiv5}.
%
%
%
\long\def\psboxit#1#2{%
\begingroup\setbox0=\hbox{#2}%
\dimen0=\ht0 \advance\dimen0 by \dp0%
    \hbox{%
    \copy0%
    }
\endgroup%
}
%
%
%
\def\Gbox#1{\psboxit{box 0.7 setgray fill}{#1}}
\def\SetTableau#1#2#3#4{%
  \gdef\Tabvrule{\vrule\vrule width-0.4pt}
  \gdef\Tabhrule{\hrule\hrule height-0.4pt}  
  \gdef\Tabstrut{\vrule height#1 depth#2 width0pt\relax}
  \gdef\Tabbox##1{\hbox to #3{\hskip0.4pt\hfill\Tabstrut$#4##1$\hfill}}
} 
%
\def\TasseTableau{\SetTableau{1.48ex}{0.32ex}{1.8ex}{\scriptstyle}}
\def\PetitTableau{\SetTableau{1.65ex}{0.55ex}{2.2ex}{\scriptstyle}}
\def\NormalTableau{\SetTableau{2.25ex}{0.75ex}{3ex}{}}
%
\def\Case#1{\vcenter{\Tabhrule%
                   \hbox{\Tabvrule\Tabbox#1\Tabvrule}\Tabhrule}}
\def\SansBord#1{\omit$\vcenter{\Tabbox#1}$}
%
\def\BordGauche#1{\omit$\vcenter{\hbox{\Tabvrule\Tabbox#1}}$}
%
\def\BordInf#1{\omit$\vcenter{\Tabbox#1\Tabhrule}$}
%
\def\CaseGrise#1{\omit\Gbox{$\Case{#1}$}}
%
\def\GenTab#1{\vcenter{\halign{&$\Case{##}$\cr#1}}\egroup}
\def\Tableau{%
  \bgroup%
  \let\ =\omit%
  \let\\=\cr%
  \offinterlineskip\GenTab}

\PetitTableau%
\newdimen\vcadre\vcadre=0.2cm 
\newdimen\hcadre\hcadre=0.2cm 
\def\GrTeXBox#1{\vbox{\vskip\vcadre\hbox{\hskip\hcadre%
      $\Tableau{#1}$\hskip\hcadre}\vskip\vcadre}}%
\def\arx#1[#2]{\ifcase#1 \relax \or%
  \ar @{-}[#2]  \or%
  \ar @{--}[#2]\fi}%

\begin{figure}[ht]
$$\xymatrix@R=0.3cm@C=0.5cm{%
 & *{\GrTeXBox{1\\}}\arx1[dd]& \\
*{\emptyset}&  & *{\GrTeXBox{2\\1\\}}& \\
 & *{\GrTeXBox{2\\}}& \\
}
\hskip3cm
\xymatrix@R=0.5cm@C=0.8cm{%
 & *{\GrTeXBox{1\\}}\arx1[d]& *{\GrTeXBox{2\\1\\}}\arx1[d]& \\
*{\emptyset}& *{\GrTeXBox{2\\}}\arx1[d]\arx2[r]& *{\GrTeXBox{3\\1\\}}\arx1[d]& *{\GrTeXBox{3\\2\\1\\}}& \\
 & *{\GrTeXBox{3\\}}& *{\GrTeXBox{3\\2\\}}& \\
}$$
  \caption{The quivers of $H_3(0)$ and $H_4(0)$.}\label{fig:quiv34}
\end{figure}
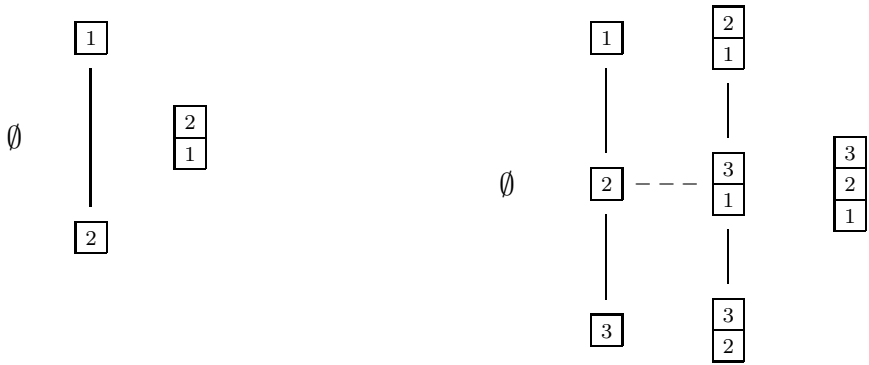


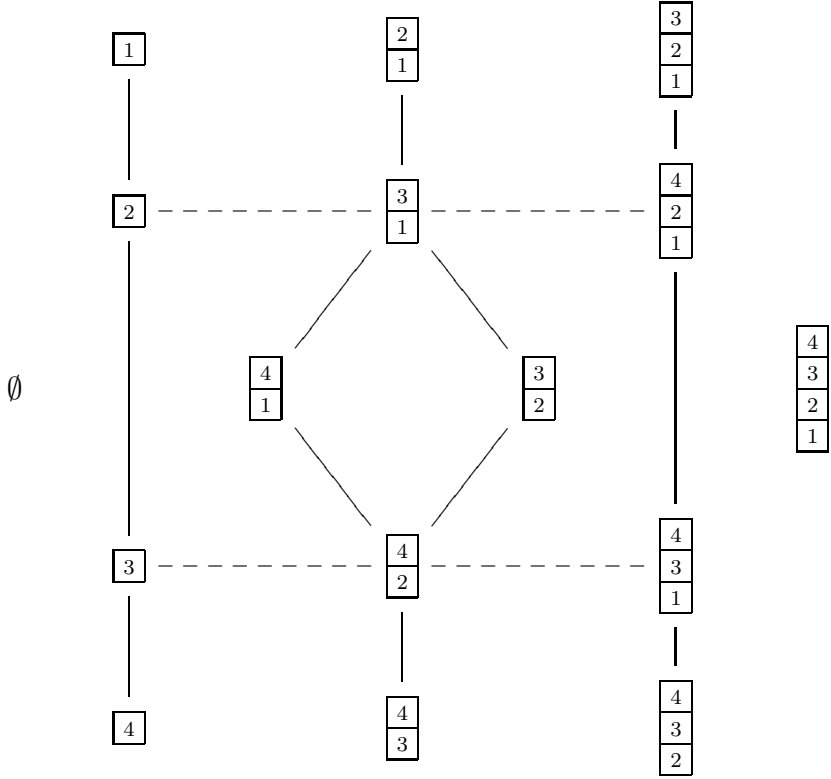
\begin{figure}[ht]
$$\xymatrix@R=0.5cm@C=1cm{%
 & *{\GrTeXBox{1\\}}\arx1[d]&  & *{\GrTeXBox{2\\1\\}}\arx1[d]&  & *{\GrTeXBox{3\\2\\1\\}}\arx1[d]& \\
 & *{\GrTeXBox{2\\}}\arx1[dd]\arx2[rr]&  & *{\GrTeXBox{3\\1\\}}\arx1[ld]\arx1[rd]&  & *{\GrTeXBox{4\\2\\1\\}}\arx1[dd]\arx2[ll]& \\
*{\emptyset}&  & *{\GrTeXBox{4\\1\\}}\arx1[rd]&  & *{\GrTeXBox{3\\2\\}}\arx1[ld]&  & *{\GrTeXBox{4\\3\\2\\1\\}}& \\
 & *{\GrTeXBox{3\\}}\arx1[d]\arx2[rr]&  & *{\GrTeXBox{4\\2\\}}\arx1[d]&  & *{\GrTeXBox{4\\3\\1\\}}\arx1[d]\arx2[ll]& \\
 & *{\GrTeXBox{4\\}}&  & *{\GrTeXBox{4\\3\\}}&  & *{\GrTeXBox{4\\3\\2\\}}& \\
}$$
  \caption{The quiver of $H_5(0)$.}\label{fig:quiv5}
\end{figure}

\noindent The vertices of the quivers are labelled by descent sets,
depicted as column shaped tableaux, instead of the corresponding
compositions. This is to emphasize the curious fact that the subgraph
on tableaux of a given height can be interpreted as the crystal graph
of a fundamental representation of $\gl_n$, or as the graph of the
Bruhat order on the Schubert cells of a Grassmannian.
 
For $n\ge 3$, $H_n(0)$ has always three blocks, a large non trivial
one, corresponding to the central connected component of the quiver,
and two one-dimensional blocks, corresponding to the two simple
projective modules $\SS_n$ and $\SS_{1^n}$. We denote by $\Gamma_n$
the quiver of the non trivial block.

For $n=2$, $H_2(0)=\C\SG_2$ is semi-simple. For $n=3$, $\Gamma_3$ is
of type $A_2$. From the well-known representation theory of such
quivers, we see that $H_3(0)$ has only 6 indecomposable modules: the 4
simple modules $\SS_I$, $I\models 4$, and the two non-simple
indecomposable projective modules $\P_{21}$ and $\P_{12}$.

For $n=4$, $\Gamma_4$ is of type $\tilde D_5$. This allows us to
conclude that $H_4(0)$ is not of finite representation type. Indeed,
choosing an orientation of $\tilde D_5$ such that the corresponding
path algebra is a quotient of $H_4(0)$, for example

\begin{center} 
\epsfig{file=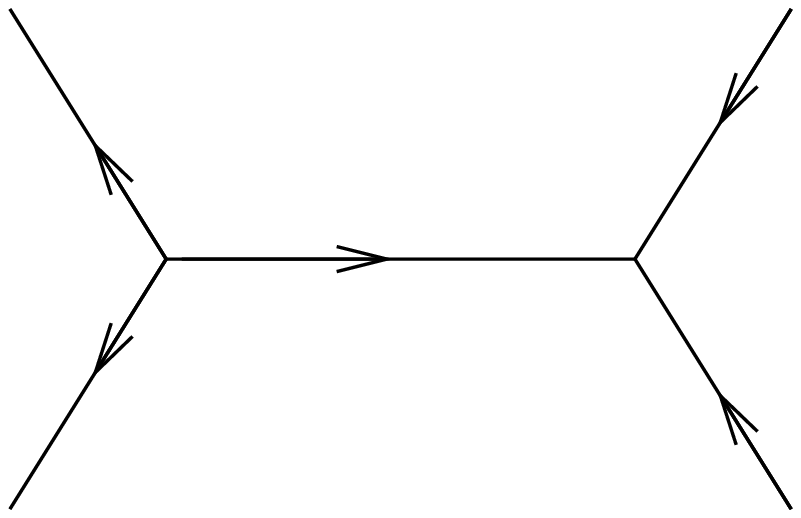}
\end{center}

\noindent (no path of length $>1$), and according to a result
of Kac \cite{Kac}, there is at least one indecomposable representation
of dimension $\alpha$ for each positive root $\alpha$, and there is an
infinite number of them.

For $n\ge 5$, then $\Gamma_n$ is considerably more complicated, and
does not belong to any familiar class of quivers. Anyway, since
$H_4(0)$ is a quotient of $H_n(0)$ for $n>4$, all these algebras are
of infinite representation type.


The quiver $\Gamma_n$ can nevertheless be described for all $n$. Indeed,
since the simple modules are one-dimensional, non trivial extensions
$$
0\rightarrow \SS_J \rightarrow M \rightarrow \SS_I \rightarrow 0
$$
are in one-to-one correspondence with indecomposable two-dimensional
modules $M$ such that $\soc M = \SS_J$ and $M/\rad M=\SS_I$.

Let $M$ be such a module, and denote by $t_i$ the matrix of $T_i$
in some basis $\{u,v\}$ of $M$. Then, $M$ is decomposable if and only if
all the $t_i$ commute. If it is not the case, let $i$ be the smallest integer
such that $t_i$ does not commute with $t_{i+1}$. The restriction of $M$
to the subalgebra $H_3(0)$ generated by $T_i$ and $T_{i+1}$ is indecomposable,
and must therefore be isomorphic to $\P_{21}$ or to $\P_{12}$. In both cases,
it is possible to choose the basis such that the matrix of $T_i$  be
\begin{equation}
  t_i=\left(\begin{matrix} -1 & 0 \\ 0 & 0 \end{matrix}\right)
\end{equation}
and $t_{i+1}$ is either
$$
\left(\begin{matrix} 0 & 0 \\ 1 & -1 \end{matrix}\right)
\quad {\rm or} \quad
\left(\begin{matrix} 0 & 1 \\ 0 & -1\end{matrix}\right)
$$
Next, $t_{i+2}$ commutes with $t_i$ and satisfies the braid relation
with $t_{i+1}$. This implies that it is either scalar (with eigenvalue
0 or $-1$) or equal to $t_{i}$. For $j>2$, $t_j$ commutes with
$t_i$ and $t_{i+1}$, so it must be a scalar matrix, again with eigenvalue
$0$ or $-1$. Also, the matrices $t_k$ for $k<i$ commute with $t_{i+1}$
and have to be scalar.

From these considerations, one obtains a complete list of indecomposable
two dimensional modules, and the following description of $\Gamma_n$

\begin{theorem}
  There is an arrow $A\rightarrow B$ between two subsets $A,B$ of
  $\{1,\ldots,n-1\}$ if and only if one of the two subsets is
  obtained from the other
  \begin{itemize}
  \item[(1)] by replacing an element $i$ either by $i+1$ or $i-1$, 
    
  \item[(2)] by deleting $i$ and inserting $i-1$ and $i+1$ if none of
    them were already present,
    
  \item[(3)] by deleting a pair $i-1$, $i+1$ and inserting $i$, if it
    was not already there. 
  \end{itemize}
\end{theorem}

From this, it is easy to see that the total number of 2-dimensional
indecomposable $H_n(0)$ modules is $(3n-7)\cdot 2^{n-3}$ for $n\ge 3$.

\subsection{Syzygies}

A way to generate infinite families of non isomorphic indecomposable
modules is to calculate the syzygies $\Omega^k\SS_I$ and
$\Omega^{-k}\SS_I$ of the simple modules $\SS_I$. The dimensions and
composition factors of these modules can be read off from the
$q$-Euler characteristics $\chi_q(\SS_I,\SS_J)$, where

\begin{equation}
  \chi_q(M,N)=\sum_{k\ge 0} q^k \dim\Ext^k_{H_n(0)}(M,N)
\end{equation}                                                                            

Indeed, if
\begin{equation}
  0\leftarrow \SS_I\leftarrow P^0\leftarrow  P^1\leftarrow P^2\leftarrow
  \cdots\leftarrow P^k\leftarrow\cdots
\end{equation}
is a minimal projective resolution of $\SS_I$, and if we write
\begin{equation}
  P^k \simeq \bigoplus \P_J^{\oplus m_{IJ}^k}
\end{equation}
then
\begin{equation}
  \sum_{k\ge 0} m_{IJ}^k q^k = \chi_q(\SS_I,\SS_J)\,.
\end{equation}
Also, since $\Omega^k\SS_I=\ker (P^{k}\rightarrow P^{k-1})$, we
have
\begin{equation}
  \ch (\Omega^k\SS_I)=\ch (P^{k-1})-\ch(\Omega^{k-1}\SS_I) \,.
\end{equation}

Moreover, for $n\le 4$, we have the more precise information
\begin{equation}
  \ch_q(\Omega^k\SS_I)= q^{-1}(\ch_q(P^{k-1})-\ch_q(\Omega^{k-1}\SS_I))
\end{equation}
where $\Omega^0\SS_I=\SS_I$ and $P^0=\P_I$. 
This formula is equivalent to the following property. Let
\begin{equation}
  A(q)=(a_{IJ}(q))_{I,J}=C(q)^{-1}\,.
\end{equation}

\begin{proposition}
  For $n\le 4$, the Poincar\'e series of
  $\Ext^*_{H_n(0)}(\SS_I,\SS_J)$ is given by
  \begin{equation}
    \chi_q(\SS_I,\SS_J)= a_{IJ}(-q)\,.
  \end{equation}
\end{proposition}
\Proof For $n=2$, this is trivial, and for $n=3$ the direct
calculation of both sides is straightforward. So let us suppose $n=4$.
As we have seen, $H_4(0)$ is a self injective algebra of infinite
representation type, with radical cube $0$, but radical square
nonzero. Hence, we can apply Theorem 1.5 of \cite{MV}, and conclude
that $H_4(0)$ is a Koszul algebra. Also, we know that the ideal $I$
such that $H_4(0)=\C Q/I$ is graded (it is generated by a set of words
of lengths 2 or 3). Then, Theorem 5.6 of \cite{GMV1} implies the
required equality. \qed

\begin{example}
  For $n=4$, the nontrivial part of the $q$-Cartan matrix,
  corresponding to the compositions $(31),(22),(211),(13),(121),(112)$
  (in this order) is
  $$
  C(q)=
  \left[
    \begin{matrix}
      1 & q & 0 & q^2 & 0 & 0 \cr
      q&1+q^2&0 & q   & q & 0 \cr
      0 & 0 & 1 & 0   & q&q^2\cr
      q^2&q&0&1&0&0\cr
      0&q&q&0&1+q^2&q\cr
      0&0&q^2&0&q&1\cr
    \end{matrix}
  \right]
  $$
  and its inverse is
  {\footnotesize
    $$
    \hss
    \frac{1}{(1-q^2)(1-q^6)}
    \left[
      \begin{matrix}
        1 & -q(1+q^4) & -q^3 & q^6 & q^2(1+q^2) & -q^3\cr
        -q(1+q^4) & (1+q^2)(1+q^4) & q^2(1+q^2) & -q(1+q^4) &-q(1+q^2)^2&q^2(1+q^2)\cr
        -q^3&q^2(1+q^2)&1&-q^3&-q(1+q^4)& q^6\cr
        q^6&-q(1+q^4)&-q^3&1&q^2(1+q^2)&-q^3\cr
        q^2(1+q^2)&-q(1+q^2)^2&-q(1+q^4)&q^2(1+q^2)&(1+q^2)(1+q^4)&-q(1+q^4)\cr
        -q^3&q^2(1+q^2)&q^6&-q^3&-q(1+q^4)&1\cr
      \end{matrix}
    \right]\hss
    $$
  }
  By taking the Taylor expansions in the first row, one can read the
  minimal projective resolution of $\SS_{31}$. The complex is
  naturally encoded by the noncommutative symmetric function
  ${\cal P}_q(\SS_{31})$:
  $$
  (1-q^2)^{-1}(1-q^6)^{-1}(R_{31}+q(1+q^4)R_{22}
  +q^3R_{211}+q^6 R_{13}+q^2(1+q^2)R_{121}+q^3R_{112})
  $$
  $$
  =R_{31}+qR_{22}+q^2(R_{31}+R_{121})+q^3(R_{22}+R_{211}+R_{112})
  +q^4(R_{31}+2R_{121}) + O(q^5)
  $$
  so that the beginning of the resolution is
  $$
  0\leftarrow \SS_{31}\leftarrow \P_{31}\leftarrow \P_{22}\leftarrow 
  \P_{31}\oplus \P_{121}\leftarrow  \P_{22}\oplus\P_{211}\oplus\P_{112}
  \leftarrow
  \P_{31}\oplus 2\P_{121}\leftarrow \cdots
  $$
  and the $q$-characteristics of the successive syzygy modules
  are
  $\ch_q(\Omega\SS_{31})=F_{22}+qF_{13}$, $\ch_q(\Omega^2\SS_{31})
  =F_{31}+F_{121}+qF_{22}$, $\ch_q(\Omega^3\SS_{31})=F_{22}+F_{211}
  +F_{112}+q(F_{13}+F_{121})$, $\ch_q(\Omega^4\SS_{31})=
  F_{31}+2F_{121}+q(F_{22}+F_{211}+F_{112})$, and so on.
  The dimensions of these modules are given by the generating
  function
  $$
  \sum_{k\ge 0}q^k \dim\Omega^k\SS_{31}=
  \frac{(1+q)(1+q^2)}{(1-q)(1-q^3)}=1+2q+3q^2+5q^3+6q^4+7q^5 +\cdots
  $$
\end{example}

\begin{example}
  The minimal projective resolution of the Specht module $V_{22}(0)$
  is encoded by the noncommutative symmetric function
  $$
  {\cal P}_q(V_{22}(0))=
  (1-q^6)^{-1}(R_{121}+qR_{112}+qR_{211}+q^2R_{121}+q^3R_{22}
  +q^4R_{13}+q^4R_{31}+q^5R_{22})
  $$
  which has period $6$, and whose commutative image tends to
  $s_{22}$ for $q\rightarrow -1$. The Poincar\'e
  series of the $\Ext$ (Yoneda) algebra
  of $V_{22}(0)$  is
  $$
  \chi_q(V_{22},V_{22})=\frac{1+q^5}{1-q^6}\,.
  $$
\end{example}

          \section{Matrix quasi-symmetric functions}

\newcommand{\hopf}{{\cal H}}      
\newcommand{\conc}{\mathrm{conc}}
\newcommand{\produ}{\times}
\newcommand{\Produ}{\mu}
\newcommand{\cprod}{\delta}
\newcommand{\Cprod}{\mathbf{c}}
\newcommand{\dtb}{\underline\cprod}
\newcommand{\cunit}{\epsilon}
\newcommand{\Cunit}{\mathbf{e}}
\newcommand{\apode}{\alpha}
\newcommand{\Apode}{\mathbf{a}}
\newcommand{\pP}{\operatorname{\cdot}}
\newcommand{\vide}{e}
\newcommand{\convol}{\operatorname{Convol}}
\newcommand{\Endom}{\operatorname{End}}
\newcommand{\lon}{\operatorname{\ell}}
\newcommand{\mm}{\operatorname{\bf m}}
\newcommand{\GL}{{GL}}
\newcommand{\caract}{\operatorname{ch}}

\newcommand{\partof}{\vdash}                    
\newcommand{\compof}{\vDash}                    

\subsection{Definition}

To define our next generalization, we start from a totally ordered set
of commutative variables $X=\{x_1<\cdots<x_n\}$ and consider the ideal
$\C[X]^+$ of polynomials without constant term.  We denote by
$\C\{X\}=T(\C[X]^+)$ its tensor algebra. The product of this algebra
will be denoted by $\Produ$. 

In the sequel, we will consider tensor products of elements of this
algebra. To avoid confusion, we denote by ``$\pP$'' the tensor product
of the tensor algebra and call it the dot product. We reserve the
notation $\tensor$ for the external tensor product. The reader should
keep in mind that, in an expression of the form $\mm = m_1\pP
m_2\cdots m_k$, none of the $m_i$ are constant monomials. Such a
product is said to be in \emph{normal form}. Otherwise we rather write
$\mm = \Produ(m_1, m_2,\cdots, m_k)$.

A natural basis of $\C\{X\}$ is formed by dot products of
monomials (called \emph{multiwords} in the sequel),
which can be represented by nonnegative integer matrices
$M=(m_{ij})$, where $m_{ij}$ is the exponent of the variable $x_i$
in the $j$th factor of the tensor product. Since constant monomials
are not allowed, such matrices have no zero column. We say that
they are \emph{horizontally packed}. 
A multiword $\mm$ can be conveniently encoded in the following
way. Let $A$ be the \emph{support} of $\mm$, that is, the set
of those variables $x_i$ such that the $i$th row of $M$ is non zero,
and let $P$ be the matrix obtained form $M$ by removing the null rows.
We set $\mm=A^P$. A matrix such as $P$, without zero rows or
columns, is said to be \emph{packed}.

For example the multiword $\mm = a\pP ab^3e^5\pP a^2d$ is
encoded by
$\indexmat\SMat{1&1&2\\0&3&0\\0&0&0\\0&0&1\\0&5&0}$.
Its support is the set $\{a,b,d,e\}$, and the associated packed matrix
is $\SMat{1&1&2\\0&3&0\\0&0&1\\0&5&0}$.
\bigskip

Let $\MQSym(X)$ be the linear subspace of $\C\{X\}$ spanned
by the elements
\begin{equation}
\MS_M = \sum_{A\in{\cal P}_k(X)}A^M
\end{equation}
where ${\cal P}_k(X)$ is the set of $k$-element subsets of $X$, and
$M$ runs over packed matrices of height $h(m)<n$.
\bigskip

For example, on the alphabet $\{a<b<c<d\}$
\renewcommand{\indexmat}%
    {\smallmatrice{\Hack a\\\Hack b\\\Hack c\\\Hack d\\}}
$$
     \MS_\SMat{1&1&2\\0&3&0\\0&0&1}=
     \indexmat\SMat{1&1&2\\0&3&0\\0&0&1\\0&0&0}+
     \indexmat\SMat{1&1&2\\0&3&0\\0&0&0\\0&0&1}+
     \indexmat\SMat{1&1&2\\0&0&0\\0&3&0\\0&0&1}+
     \indexmat\SMat{0&0&0\\1&1&2\\0&3&0\\0&0&1}
$$
\bigskip

\begin{proposition} $\MQSym$ is a subalgebra of $\C\{X\}$. Actually,
$$
\MS_P\MS_Q =\sum_{R\in\ash (P,Q)} \MS_R
$$
where the {\em augmented shuffle} of $P$ and $Q$,
$\ash(P,Q)$ is defined as follows: let $r$  be an
integer between $\max(p,q)$ and $p+q$, where $p=h(P)$ and $q=h(Q)$.
Insert null rows in the matrices $P$ and $Q$ so as to form 
matrices $\tilde P$ and $\tilde Q$ of height $r$. Let $R$ be
the matrix $(\tilde P,\tilde Q)$. The set $\ash (P,Q)$ is formed
by all the matrices without null rows obtained in this way.
\end{proposition}
For example :
$$
\begin{array}{l}
  \MS{\SMat{2&1\\1&0}}\MS_{\SMat{3&1}} = \\[3mm]
\qquad\ 
  \MS{\SMat{2&1&0&0\\1&0&0&0\\0&0&3&1}}+\MS{\SMat{2&1&0&0\\1&0&3&1}}+
  \MS{\SMat{2&1&0&0\\0&0&3&1\\1&0&0&0}}+\MS{\SMat{2&1&3&1\\1&0&0&0}}+
  \MS{\SMat{0&0&3&1\\2&1&0&0\\1&0&0&0}}
\end{array}
$$
\pagebreak[3]

Let us endow $\MQSym$ with a Hopf algebra structure.  Let
$Y=\{y_1<\cdots<y_n\}$ be a second totally ordered set of variables,
of the same cardinality as $X$. We identify the tensor product
$\MQSym(X)\otimes \MQSym(X)$ with $\MQSym(X\oplus Y)$, where $X\oplus
Y$ denotes the ordered sum of $X$ and $Y$. The natural embedding
\begin{equation}
\Delta\ :\ \MQSym(X)\longrightarrow
\MQSym(X\oplus Y) \simeq \MQSym(X)\otimes \MQSym(X)
\end{equation}
defined by $\Delta(\MS_M(X))=\MS_M(X\oplus Y)$ can be interpreted
as a comultiplication.
\medskip

For example
$$
\begin{array}{r@{\,}l}
\Delta\left(\MS{\SMat{1&0&3\\ 0&2&1\\ 0&0&3\\ 1&0&2}}\right) = 
&1\otimes \MS{\SMat{1&0&3\\ 0&2&1\\ 0&0&3\\ 1&0&2}} 
 \ +\ \MS_{\SMat{1&3}}\otimes\MS{\SMat{0&2&1\\ 0&0&3\\ 1&0&2}} 
 \ +\ \MS{\SMat{1&0&3\\ 0&2&1}}\otimes\MS{\SMat{0&3\\ 1&2}} \\
&\ +\ \MS{\SMat{1&0&3\\ 0&2&1\\ 0&0&3}}\otimes\MS_{\SMat{1&2}} 
 \ +\ \MS{\SMat{1&0&3\\ 0&2&1\\ 0&0&3\\ 1&0&2}}\otimes 1 
\end{array}
$$

\emph{From now on, unless otherwise stated, we suppose that $X$ is
infinite.}

Let $\mu\ :\ f\otimes g \mapsto fg$ be the multiplication
of $\MQSym$ (induced by the multiplication of the tensor algebra), and
let $e$ be the restriction to $\MQSym$ of the augmentation of
$T(\C[X]^+)$. Introduce a grading by setting $\deg(\MS_M)=\sum m_{ij}$
and denote by $\MQSym_d$ the homogeneous component of degree $d$.

\begin{proposition}
$(\MQSym,\mu,1,\Delta,e)$ is a self dual graded bialgebra, 
the duality pairing being given by $\pairing{\MS_P}{\MS_Q}=\delta_{P,{}^tQ}$.
\end{proposition}

The Hilbert series of $\MQSym$ can be expressed directly or in terms
of scalar products of ribbon Schur functions. One has
\begin{equation}
\dim(\MQSym_d)
   = \sum_{l>0,\ h>0} \binom{d+lh-1}{lh-1}
      \left(\frac{1}{2}\right)^{l+h+2}
   =\sum_{|I|=|J|=d} 2^{2d-\ell(I)-\ell(J)}\pairing{r_I}{r_J}
\end{equation}
which yields 
$$
\begin{array}{r@{\,}l}
\sum_{d\ge 0} &\dim(\MQSym_d)\, t^d \\
&=1+t+5t^2+33t^3+281t^4+2961t^5+37277t^6+546193t^7+9132865 t^8+\cdots
\end{array}
$$

To a packed matrix $P$, we can associate two compositions $I =
\text{Row}(P)$ and $J = \text{Col}(P)$ formed by the row-sums and
column-sums of $P$. The Hilbert series in an easy consequence of the
classical fact that packed matrices of degree $d$ are in bijection
with double cosets of $\SG_J\backslash\SG_d/\SG_I$ where
$\SG_{(i_1,i_2,\dots,i_k)}$ is the Young subgroup
$\SG_{i_1}\times\SG_{i_2}\times\dots\times\SG_{i_k}$. It is well
known (\cf \cite{Ge}) that their number is
$2^{2d-\ell(I)-\ell(J)}\pairing{r_I}{r_J}$. 
\bigskip

Let ${\rm Ev}$ be the linear map defined by 
\begin{equation}
\begin{array}{rccl}
  {\rm Ev} \ : & \MQSym & \longrightarrow & QSym \\
               & \MS_P & \longmapsto & {\rm Ev}(\MS_P)=M_{{\rm Row}(P)}
\end{array}
\end{equation}
\begin{proposition}
${\rm Ev}$ is an epimorphism of bialgebras. Dually,
the transposed map 
\begin{equation}
\begin{array}{rccl}
  {}^t{\rm Ev}\ :  & \Sym &\rightarrow &\MQSym \\
                   & \MS_P & \longmapsto &
                           {}^t{\rm Ev}(S^I)=\sum_{{\rm Col}(P)=I}\MS_P 
\end{array}
\end{equation}
is a monomorphism of bialgebras.
\end{proposition}

Therefore, $\MQSym$ admits $QSym$ as a quotient and $\Sym$ as a subalgebra.
The basis $\MS_P$ can be regarded as a simultaneous generalization
of the dual bases $M_I$ and $S^I$.
Moreover, $\C\free{X}$ is naturally a subalgebra of $\C\{X\}$,
words being identified with multiwords with exponent matrix
having only one $1$ in each column. It is clear that this embedding maps
$\FQSym$ to a subspace of $\MQSym$.

\subsection{Algebraic structure}

We now elucidate  the structure of $\MQSym$ as an algebra. 
To describe a generating family, we need the following definitions.
Let $P$ be a packed matrix of height $h$. To a composition
$K=(k_1,\ldots,k_p)$ of $h$, we associate the matrix $P\rangle K $
defined as follows. Let $R_1,R_2,\ldots,R_h$ be the rows of $P$.
The first row of $P\rangle K$ is the sum of the first $k_1$ rows of $P$,
the second row is the sum of the next $k_2$ rows of $P$, and so on.
We end therefore with a matrix of height $p$. For example,
$$
\left.\SMat{1&2&0&2\\0&1&2&1\\1&2&0&0\\0&3&1&5\\1&3&1&0}\right\rangle
   (3,2)=\SMat{2&5&2&3\\1&6&2&5}\,.
$$
Generalizing the idea of \cite{MR1}, we set
\begin{equation}
\FS_P =\sum_{|K|=h} \frac{1}{K!} \MS_{P\> K}
\end{equation}
where $K!=k_1!k_2!\cdots k_p!$. The family $\{\FS_P\}$, where
$P$ runs over packed matrices, is a homogeneous basis of $\MQSym$.

Let us say that a packed matrix $A$ is \emph{connected} if it cannot
be written in block diagonal form 
$$
A=\left(\begin{matrix}B&0\cr 0&C\cr\end{matrix}\right)
$$
where $B$ and $C$ are not necessarily square matrices.

\begin{theorem}
$\MQSym$ is freely generated by the family $\{\FS_A\}$, where $A$ runs
over the set of connected packed matrices.
\end{theorem}

\subsection{Convolution}

The goal of this subsection is to find an interpretation of $\MQSym$
in terms of invariant theory. The first part of the construction
applies to any Hopf algebra.

\subsubsection*{Hopf Algebra background}

First, let $(\hopf,1,\Produ,\cprod,\cunit,\apode)$ be a graded Hopf
algebra. One can define a bialgebra structure on the augmentation ideal
$T(\hopf^+)$. The coproduct is $\Cprod$ is defined as follows. Let
$\mm=m_1\pP m_2\cdots m_p$ be a normal form dot product ($m_i\in
\hopf^+$). Let $\cprod(m_i)=\sum m'_i\tensor m''_i$. Then, one sets
\begin{equation}\label{DEF_CoprodTens}
  \Cprod(\mm)=\sum
     \Produ\left( m'_1, m'_2,\cdots, m'_p \right)
        \tensor
     \Produ\left(m''_1, m''_2,\cdots, m''_p\right)\,.
\end{equation}
For example, with $\hopf = \C[x,y]$, one has
$$
\Cprod(x^2\pP y)= x^2\pP
y\tensor 1 + 2 x\pP y\tensor x + y\tensor x^2 + x^2\tensor y + 2
x\tensor x\pP y + 1\tensor x^2\pP y\,.
$$
If $\cprod$ is cocommutative, then $\Cprod$ is obviously so.
The co-unit is the coordinate of the empty tensor
\begin{equation}
  \Cunit(1)=1
  \qquad\text{and}\qquad
  \Cunit(m_1\pP m_2\cdots m_p)=0\text{ if }p>0\text{ and }m_i\in\hopf^+\,.
\end{equation}
If $\hopf$ is graded, one defines a gradation on $T(\hopf^+)$ by
\begin{equation}
   \deg(m_1\pP m_2\cdots m_p)=\deg(m_1)+\deg(m_2)+\cdots +\deg(m_p)\,,
\end{equation}
Note that this gradation differs from the standard one on
tensors, which we will call length
\begin{equation}
   \lon(m_1\pP m_2\cdots m_p)=p\,. 
\end{equation}
\bigskip

Now, $f$ and $g$ being two endomorphisms of $\hopf$, the convolution
of $f$ and $g$ is defined by $f*g=\Produ\circ(f\tensor g)\circ\cprod$.
In the sequel, all the endomorphisms will be homogeneous. The
convolution of an endomorphism of degree $p$ with an endomorphism of
degree $q$ is of degree $p+q$. One denotes by $\convol(\hopf)$ the
\emph{convolution algebra of the homogeneous endomorphisms of $\hopf$}
and by $\Endom^n(\hopf)$ the vector space of \emph{homogeneous
  endomorphisms of degree~$n$}.

\subsubsection*{Operator associated with a packed matrix}

To each packed matrix $A$ of total sum $n$ we associate a canonical
endomorphism $f_A$ of the $n^{\text{th}}$ homogeneous component of
$T(\hopf^+)$. First of all, let $K=(k_1,\ldots,k_q)\in\NN^q$. Let us
define
\begin{alignat}{2}
  \notag
\cprod^{(K)} \ :\ &&\hopf & \longrightarrow \hopf^{\tensor q} \\
                  &&  m   & \longmapsto
    (\pi_{k_1}\tensor\cdots\tensor\pi_{k_q})\circ\cprod^q(m)
\end{alignat}
where $\pi_d$ is the projector on the homogeneous component of degree
$d$ of $\hopf$. Thus $\cprod^{(K)}$ takes an element of degree $|K|$
and sends it to an element of $\hopf^{\tensor q}$ of degree
$(k_1,\ldots,k_q)$, killing all components of other degrees.

Let $A=(a_{i,j})$ be a packed $p\times q$ matrix of total sum $n$. The
row sum of $A$ is a composition $r=(r_1,\ldots,r_p)$ and the column
sum is $c=(c_1,\ldots,c_q)$. Let us denote by $R_1,\dots,R_p$ the rows
of $A$. Finally suppose that $\mm = m_1\pP m_2\cdots m_r$ is an
element of $T(\hopf^+)$. Then we define $f_A$ by
\begin{equation}
  f_A(m_1\pP m_2\cdots m_r)=
\begin{cases}
   \Produ_q^p(\cprod^{L_1}(m_1),\cdots,\cprod^{L_p}(m_p))
     & \text{if $r=p$,} \\
   0 & \text{otherwise,}\\
\end{cases}
\end{equation}
where $\Produ_q^p=(\Produ^p)^{\tensor q}$ is the product of the $p$
tensor $\cprod^{L_i}(m_i)$ of length $q$. Thus we get an element of
$(\hopf^{+})^{\tensor q}$ of degree $(c_1,\ldots,c_q)$. Remark
that $f_A(\mm)$ is null unless $m$ is of degree $l=(l_1,\ldots,l_p)$. 

\begin{example}
  let $A=\SMat{2&0&1\\0&2&3}$. The associated morphism $f_A$ kills all
  tensors of degree different from $(3,5)$. Let $\mm=abc \pP a^4b$. 
  Then
\begin{eqnarray*}
  \cprod^{(2,0,1)}(abc)&=&
    ab\tensor1\tensor c+
    ac\tensor1\tensor b+
    bc\tensor1\tensor a,\\
  \cprod^{(0,2,3)}(a^4b)&=&
    \binom{4}{2}(1\tensor a^2\tensor a^2b)+
    \binom{4}{1}(1\tensor ab\tensor a^3)\,.
\end{eqnarray*}
Finally,
\begin{equation*}
  \begin{array}{c@{\,}r@{\,}c@{\,}c@{\,}c@{\,}c@{\,}c@{\,}l}
  F_A(abc\pP a^4b) &=
     6\,(&ab\pP a^2\pP a^2bc &+& ac\pP a^2\pP a^2b^2 &+& bc\pP a^2\pP a^3b&)\\
   &+4\,(&ab\pP ab \pP a^3c &+& ac\pP ab \pP a^3b   &+& bc\pP ab \pP a^4&)\,.
  \end{array}
\end{equation*}
\end{example}

The following example, which is some sense generic, is of crucial
importance.
\newcommand{\KK}{{\mathbb K}}                   
\newcommand{\X}{X}
\newcommand{\card}{\operatorname{\#}}           
\newcommand{\mat}[1]{\{#1\}}                    
\newcommand{\pol}[1]{[#1]}                      

\begin{example}\label{EXPL_GenMult}
  
  Let $A=(a_{i,j})$ of size $p\times q$ and degree (total sum) $n$.
  Let us consider
  \begin{equation}
    \KK\mat\X=T(\KK\pol{x_1,x_2,\dots,x_n}^{+})\,.
  \end{equation}
  To the composition $r=(r_1,\ldots,r_p)$ of the row sum of $A$, we
  associate the generic multiword of degree $r$ denoted by $\mm_{(r)}$
  and defined as follows: let $d_1=r_1$, \ $d_2=r_1+r_2$, \ \dots, \ 
  $d_i=r_1+\dots+r_i$, \ \dots, and $d_p=n$ the descents of $r$. Define
  \begin{equation}
    \mm_{(r)} = 
    \left( \prod_{i=1}^{d_1} x_i \right)
      \pP\,
    \left( \prod_{i=d_1+1}^{d_2} x_i \right)
      \cdots
    \left( \prod_{i=d_{p-1}+1}^{d_{p}}x_i \right)\,,
  \end{equation}
  or, equivalently,
  \begin{equation}
    \mm_{(r)} = 
    x_1\diamond_1 x_2\diamond_2 x_3\diamond_3\cdots\ \diamond_{n-1}x_n,
  \end{equation}
  where $\diamond_i$ is the commutative multiplication if $i$ is not a
  descent of $r$, and the dot product otherwise. Let
  $\X_D=\prod_{i\in D} x_i$ where $D$ is a subset of $\{1\dots n\}$ and
  moreover, let $D_i$ denote the integer interval
  $\{d_{i-1}+1,\dots,d_i\}$. Then
  \begin{equation}
    \mm_{(r)} = \X_{D_1} \pP\, \X_{D_2} \cdots \X_{D_p}\,.
  \end{equation}
  \medskip
  
  Let us compute the image of $\mm_{(r)}$ by $f_A$. For all
  $K=(k_1,\ldots,k_p)\in\NN^q$ of sum $s$ one has
  \begin{equation}
    \cprod^{(K)}(\X_{\{u,\dots, u+s\}})
    = \sum_{I_1,\ldots,I_q}
    ( \X_{I_1}\tensor\cdots\tensor \X_{I_q} )\,,
  \end{equation}
  where the sum is over all set-partitions $I_1,\ldots,I_q$ of the
  integer interval $\{u,\dots, u+s\}$ such that $\card(I_1)=k_1$,
  \dots , $\card(I_q)=k_q$. It follows that
  \begin{equation}
    \newcommand{\ds}{\displaystyle}
    f_A(\mm_{(r)}) = \sum_{(I_{i,j})}
    \left(
      \X_{\ds\cup I_{i,1}}\tensor\cdots\tensor \X_{\ds\cup I_{i,q}}
    \right)\,,
  \end{equation}
  the sum is over all $p\times q$-matrices $(I_{i,j})$ whose entries
  are subsets of $\{1,\dots n\}$ and such that
  \begin{itemize}
  \item for all $i,j$, one has $\card(I_{i,j})=a_{i,j}$,
  \item for all $i$ the set $\{I_{i,1},\ldots,I_{i,q}\}$ defines
    a partition of the interval $D_i=[d_{i-1}+1,\dots,d_i]$.
  \end{itemize}
  For example, with $A=\SMat{0&1&1\\1&0&2}$, one has $r=(2,3)$. Then
  the generic multiword $\mm_{(r)}$ reads
  \begin{equation}
    \mm_{(2,3)}=x_1x_2\pP x_3x_4x_5\,.
  \end{equation}
  Then 
  \begin{equation*}
    \begin{array}{crcccccl}
      \cprod^{(0,1,1)}(x_1x_2)  &=&
      1\tensor x_1\tensor x_2 &+& 
      1\tensor x_2\tensor x_1,\\
      \cprod^{(1,0,2)}(x_3x_4x_5)  &=&
      x_3\tensor 1\tensor x_4x_5 &+& 
      x_4\tensor 1\tensor x_3x_5 &+& 
      x_5\tensor 1\tensor x_3x_4\,.
    \end{array}
  \end{equation*}
  Finally
  \begin{equation*}
    \begin{array}{crccccc@{}l}
      F_A(x_1x_2\pP x_3x_4x_5)
      &=& x_3\pP x_1\pP x_2x_4x_5 &+&  x_4\pP x_1\pP x_2x_3x_5 &+&
          x_5\pP x_1\pP x_2x_3x_4 \\
      &+& x_3\pP x_2\pP x_1x_4x_5 &+&  x_4\pP x_2\pP x_1x_3x_5 &+&
          x_5\pP x_2\pP x_1x_3x_4 & \,.
    \end{array}
  \end{equation*}
  
\end{example}

\begin{theorem} The map 
  \begin{alignat}{1}
     \notag
        \MQSym & \longrightarrow \convol(T(\hopf^{+})) \\
         \MS_A& \longmapsto   f_A
  \end{alignat}
  is a homomorphism of algebras.  
\end{theorem}

\Proof
\newcommand{\detass}[2]{\mathchoice{(#1\!\updownarrow\!_#2)}%
   {(#1\!\updownarrow\!_#2)}{(#1\updownarrow_#2)}{(#1\updownarrow_#2)}}
The first step of the proof is to see that the definition of the
morphism $f_A$ can be extended to non-packed matrices. If $B$ is an
integer $p\times q$ matrix the preceding definition gives a
morphism
\begin{equation}
   \tilde{f}_B\ :\ \hopf^p\longrightarrow T(\hopf^{+})\,,
\end{equation}
With this notation, one has the following easy lemma:
\begin{lemma}
   Let $A$ be a packed matrix of height $h$. Let $$\mm=\Produ(m_1,
   m_2,\dots, m_p)$$ a tensor, not necessarily in normal form ($m_i$
   can be constant). Then,
   \begin{equation}
     f_A(\mm)=\sum_{B\in \detass{A}p}\tilde{f}_{B}(m_1, m_2,\dots, m_p)\,,
   \end{equation}
   where ${B}$ runs over the set $\detass{A}p$ of matrices of height
   $p$ obtained by inserting $0$ rows in the matrix $A$. 
\end{lemma}
Now, let $\mm=m_1\pP m_2\cdots m_p$. Suppose $\cprod(m_i)=\sum m'_i\tensor
m''_i$. By definition
\begin{equation}
  \Cprod(\mm)  =  \sum
     \Produ\left( m'_1, m'_2, \cdots, m'_p \right)
  \tensor
     \Produ\left( m''_1, m''_2, \cdots, m''_p \right)\,.
\end{equation}
And therefore, if $A$ and $A'$ are two packed matrices,
\begin{equation}
  (f_{A}\tensor f_{A'})\circ\Cprod(\mm)=
    \sum_{B\in\detass{A}p,\ B'\in\detass{A'}p}
       \tilde{f}_{B} \left( m'_1, \dots, m'_p \right)
         \tensor
       \tilde{f}_{B'}\left( m''_1, \dots, m''_p \right)\,,
\end{equation}
which gives
\begin{equation}
  \Produ\circ(f_{A}\tensor f_{A'})\circ\Cprod(\mm)=
    \sum_{B\in\detass{A}p,\ B'\in\detass{A'}p}
       \tilde{f}_{BB'} \left(
         m_1, \dots, m_p
       \right)\,,
\end{equation}
where $BB'$ is the concatenation of $B$ and $B'$. This is exactly the
set of unpackings $\detass{C}p$ of the matrices $C$ appearing in the
product $\MS_A\MS_{A'}$ which acts non-trivialy on $\mm$. \cqfd

Note that the theorem is true even if $\hopf$ is not cocommutative. 

\subsubsection*{Interpretation}

First, we reformulate the definition of $\C\{X\}$ (with
$X=\{x_1,\ldots,x_n\}$) in a slightly more abstract way. Let $V$ be an
$n$-dimensional vector space with basis $X$. The polynomials in $X$ can
be seen as the symmetric algebra of $V$. The graded bialgebra
structure on $\C\{X\}=T(\C[X]^+)$ gives a structure on $T(S^+(V))$.
Moreover, since the definition of the operations does not depend on
the basis, this structure is canonical. Let $\rho$ be the natural
representation of $GL(V)$ in $\Endom(T(S^+(V)))$.

\begin{theorem} There exists a canonical homomorphism
  \begin{equation*}
    \phi\ :\ \MQSym
    \longrightarrow \Endom (T(S^+(V)))
  \end{equation*}
  from $\MQSym$ to $\Endom (T(S^+(V)))$ regarded as a convolution\
  algebra, such that for all $d$, $\phi(\MQSym_d)$
  is the commutant $\Endom_{GL(V)}(T(S^+(V))_d)$ of $\rho(GL(V))$ in
  the homogeneous component of degree $d$ of $\Endom(T(S^+(V)))$.
  Moreover, $\phi$ is one-to-one for $d\le n$.
\end{theorem}

\Proof The endomorphism $f_A$ associated with a matrix $A$ is defined
by means of the product, coproduct, and the homogeneous projector of
$T(S^+(V))$. But all these operations commute with the action of
$GL(U)$. Then $f_A$ commutes with $GL(U)$. \medskip

We will prove the theorem in two steps:
\begin{itemize}
\item In the first step, we suppose that the dimension $N$ of $V$ is
  greater than $n$. We will prove that $\phi$ is one-to-one and, by an
  argument of dimension, we get that $\phi(\MQSym_d)$ is exactly the
  commutant $\Endom_{GL(V)}(T(S^+(V))_d)$.
\item Then by a restriction argument we will conclude in every case. 
\end{itemize}
\bigskip

Let us choose a basis $X=\{x_1,\dots x_n\}$ of $V$. In
example~\ref{EXPL_GenMult}, we have computed the image of the generic
multiword $\mm_{(r)}$ by $f_A$ where $A$ is a matrix of row sum $r$.
Notice that $n$ is sufficient to express $\mm_{(r)}$ since $N$ is
bigger than the sum $n$ of $A$.

Let us recall some notation: $d_1,\dots,d_p$ denote the descents of
$r$ and $D_i$ the integer interval $[d_{i-1}+1,\dots,d_i]$. Let us
suppose that $\mm'$ is a multiword of the form
\begin{equation}
   \mm'=X_{I_1}\pP X_{I_2}\cdots X_{I_q}.
\end{equation}
where $I_1,\dots,I_q$ is a partition of the set $\{1,\dots,n\}$. There
exists only one matrix $A$ such that $\mm'$ appears in the image of
$\mm_{(r)}$ by $f_A$:
\begin{equation}
  A=
  \begin{bmatrix}
     \card(I_1\cap D_1) & \cdots & \card(I_q\cap D_1) \\
          \vdots           & \ddots &      \vdots     \\ 
     \card(I_1\cap D_p) & \cdots & \card(I_q\cap D_p) 
  \end{bmatrix}.
\end{equation}
This proves the injectivity of $\phi$.

Now we will show that the dimension of $\MQSym^n$ and of the 
commutant of $\GL(U)$ in $\Endom^n(T(S^{+}(U)))$ are equal. Let us
compute the graded character of the representation $\GL(U)$ on
$T(S^{+}(U))$. It is well known that the character of $S^d(U)$ is the
Schur function $s_{(d)}$ which is equal to the complete function
$h_{d}$.

The graded character of $S^{+}(U)$ is then:
\begin{equation}
   \caract_t(S^{+}(U))=\sum_{d>0} h_{d} t^d\,.
\end{equation}
Therefore
\begin{equation}
   \caract_t(T(S^{+}(U)))=\sum_I h_I t^{|I|}\,,
\end{equation}
where $I$ runs over the set of all compositions. Note that $h_I$
only depends on the partition associated with $I$. Then one uses the
classical identity
\begin{equation}
\sum_I h_I u^{\lon(I)}=\sum_J r_J u^{\lon(J)}(1+u)^{|J|-\lon(J)}.
\end{equation}
Now, we extract the homogeneous components, with $u=1$. This gives
\begin{equation}
   \caract_t(T(S^{+}(U)))=\sum_{J} r_J 2^{d-\lon(J)} t^{|J|} 
\end{equation}
The multiplicity of the irreducible representation $\chi_\lambda$ of
$\GL_N$ in the homogeneous component of degree $d$ of $T(S^{+}(U))$ is
therefore given by the scalar product
\begin{equation}
\sum_{\lambda\partof n}
   \sum_{I\compof n,\ J\compof n}
      2^{2n-\lon(I)-\lon(J)}
         \pairing{r_I}{s_\lambda}\pairing{s_\lambda}{r_J}.
\end{equation}
The sum is extended to partitions all $\lambda$ of length smaller than
the dimension $N$ of $U$. Thus if $N\geq n$, all the Schur functions
$s_\lambda$ appear. Moreover, since they form an orthonormal basis of
$Sym$
\begin{equation}
\sum_{\lambda\partof n}
   \pairing{r_I}{s_\lambda}\pairing{s_\lambda}{r_J}=\pairing{r_I}{r_J}.
\end{equation}
This proves that the dimension of the commutant of $\GL(U)$ in
the $n$-homogeneous space of $T(S^{+}(U))$ is equal to the dimension
of $\MQSym_n$, which implies the first part of the theorem. 
\bigskip

Now we are in the case where the dimension $N$ of the vector space $U$
is less than $n$. Let $V=U\oplus W$ be of dimension $n$.

\begin{lemma}
  Let $U\subset V$ be two vector spaces. Then the restriction
\begin{alignat*}{2}
  \text{Rest}_{U\subset V}\ :\
     && \Endom_{GL(V)}(T(S^+(V))_n)  &
        \longrightarrow \Endom_{GL(U)}(T(S^+(U))_n)\\
     && f &    \longmapsto f_{T(S^{+}(U))_n} 
\end{alignat*}
is surjective. 
\end{lemma}

Using this lemma one has that every element of
$\Endom_{GL(U)}(T(S^+(U))_n)$ is the restriction of some element of
$\rho_V(\MQSym)$. But it is clear that the endomorphism $F^U_A$
associated with a matrix $A$ on $U$ is the restriction to
$T(S^{+}(U))$ of $F^V_A$ associated with $A$ on $V$. This concludes
the proof of the theorem.  \bigskip

\newcommand{\comm}{\operatorname{Comm}}
\newcommand{\Vect}{\operatorname{Vect}}
It remains to prove the lemma. First we have to prove that the image
of an element of $T(S^{+}(U))$ by $f$ is still in $T(S^{+}(U))$.

Let us set $\comm_U = \Endom_{GL(U)}(T(S^+(U)))$.  Let $\mm$ be an
element of $T(S^{+}(V))$. Let $\Vect(\mm)$ be the smallest subspace
$W\subset V$ such that $\mm\in T(S^{+}(W))$. Then clearly if
$g\in\GL(V)$ then $\Vect(g(\mm)) = g(\Vect(\mm))$. But if $f$ commutes
with $\GL_V$, it also commutes with the projectors on $\Vect(\mm)$, so
that
\begin{equation}
  \text{for all }f\in\comm_V,
     \quad f(\mm)\in T(S^{+}(\Vect(\mm)))\,.
\end{equation}

Now, let us prove the surjectivity. Let $g\in\comm_U$. Let $\mm\in
T(S^+(V))$. Under the assumption $\dim(\Vect(\mm)) \geq \dim(u)$, one
can define the image $f(\mm)$ by conjugation as follows: choose an
injective morphism $h_{\mm}\,:\,\Vect(\mm)\mapsto U$. Then obviously
$h_{\mm}(\Vect(\mm))=\Vect(h_{\mm}(\mm))$ and one can set
\begin{equation}
f(\mm)=
  \begin{cases}
    h_{\mm}^{-1}\circ g \circ h_{\mm}(\mm) & \text{if $\dim(V(\mm))\leq u$}, \\
    0                            & \text{otherwise.} \\
  \end{cases}
\end{equation}
Since $g$ commutes with $\GL(U)$, the vector $f(\mm)=h_{\mm}^{-1}\circ
g \circ h_{\mm}(\mm)$ does not depend on the choice of $h_{\mm}$.
Hence if $\mm\in T(S^+(U))$, one can take $h_{\mm}=id_u$ and thus
$f_U=g$. Moreover it is easy to see that $f$ commutes with $\GL(V)$.
\cqfd


\bigskip
\begin{flushleft}
\noindent
\sc
Isaac Newton Institute for Mathematical Science\\
20 Clarkson road\\
Cambridge, CB3 0EH, U.K.\\
\medskip

\sf Permanent addresses :\\ 

\medskip\noindent
\sc G. Duchamp:\\ 
Laboratoire d'Informatique Fondamentale et Appliqu\'ee de Rouen\\
Facult\'e des Sciences, Universit\'e de Rouen,\\
76821 Mont-Saint-Aignan cedex,\\
FRANCE.\\

\medskip\noindent
F. Hivert and J.-Y. Thibon:\\
Institut Gaspard Monge, \\
Universit\'e de Marne-la-Vall\'ee,\\
77454 Marne-la-Vall\'ee cedex,\\
FRANCE
\end{flushleft}

\footnotesize
\begin{thebibliography}{99}
%
\bibitem{Ab}{\sc L.~Abrams}, {\it Modules, comodules and cotensor
    products over Frobenius algebras}, {\tt math.RA/9806044}.
%
\bibitem{ARS}{\sc M.~Auslander, I.~Reiten} and {\sc S.-O~Smal\o}, {\it
    Representation theory of Artin algebras}, Cambridge Studies in
  Advanced Mathematics, {\bf 36} (1997), Cambridge University Press,
  Cambridge.
%
\bibitem{Comtet}{\sc L.~Comtet}, 
{\it Sur les coefficients de l'inverse de la s\'erie formelle $\sum n!
  t^n$},
C.R. Acad. Sci. Paris {\bf A 275} (1972), 569--572.
%
\bibitem{CR}{\sc C.W.~Curtis} and {\sc I.~Reiner}, 
{\it Methods of representation theory}, 
2 vols., Wiley-Interscience, New-York, 1981.
%
\bibitem{Du}{\sc G. Duchamp}, 
{\it Orthogonal projection onto the free Lie algebra}, 
Theoret. Comp. Sci. {\bf 79} (1991), 227--239.
%
\bibitem{NCSF3} {\sc G.~Duchamp, A.~Klyachko, D.~Krob} and {\sc
    J.-Y.~Thibon}, 
{\it Noncommutative symmetric functions III{\,}:
    Deformations of Cauchy and convolution algebras},
Disc. Math. and Theor. Comput. Sci.  {\bf 1} (1997), 159--216.
%
\bibitem{DKLT}{\sc G.~Duchamp , D.~Krob, B.~Leclerc} and {\sc J.-Y.~Thibon},
{\it Fonctions quasi-sym\'etriques, fonctions sym\'etriques non
commutatives, et alg\`ebres de Hecke \`a $q=0$}, 
C.R. Acad. Sci. Paris, {\bf 322} (1996),  107--112.
%
\bibitem{FS}{\sc D.~Foata} and {\sc M.P.~Sch\"utzenberger}, 
{\it Major index and inversion number of permutations}, 
Math. Nachr. {\bf 83} (1978), 143--159.
%
\bibitem{NCSF1} {\sc  I.M.~Gelfand, D.~Krob, B.~Leclerc, A.~Lascoux,
V.S.~Retakh} and {\sc J.-Y.~Thibon}, 
{\it Noncommutative symmetric functions},
Adv. in Math., {\bf 112} (1995), 218--348.
%
%
\bibitem{Ge} {\sc I.~Gessel}, {\it Multipartite P-partitions and inner
products of skew Schur functions}, [in ``Combinatorics and algebra",
C. Greene, Ed.], Contemporary Mathematics, {\bf 34} (1984), 289--301.
%
\bibitem{GMV1}{\sc E.L. Green} and {\sc R. Mart\'\i nez-Villa},
{\it Koszul and Yoneda algebras}, in Representation Theory of Algebras,
CMS Conference Procceedings, vol. 18 (1996), 247--297.
%
\bibitem{GMV2}{\sc E.L. Green} and {\sc R. Mart\'\i nez-Villa},
{\it Koszul and Yoneda algebras II}, in Algebras and Modules II,
CMS Conference Procceedings, vol. 24 (1998), 227--244.
%
\bibitem{Kac}{\sc V.G. Kac}, 
{\it Some remarks on representations of quivers and infinite root systems}, 
Representation theory, II 
(Proc. Second Internat. Conf., Carleton Univ., Ottawa, Ont., 1979), 
pp. 311--327, Lecture Notes in Math., 832, Springer, Berlin, 1980. 
%
%
\bibitem{Kly}{\sc A.A. Klyachko}, {\it Private communication} (1995).
%
\bibitem{NCSF2} {\sc D.~Krob, B.~Leclerc} and {\sc J.-Y.~Thibon}, 
{\it Noncommutative symmetric functions II{\,}: Transformations of
alphabets}, 
Int. J. of Alg. and Comput. {\bf 7} (1997), 181--264.
%
\bibitem{NCSF4} {\sc D.~Krob} and {\sc J.-Y.Thibon},
{\it Noncommutative symmetric functions IV{\,}: 
  Quantum linear groups and Hecke algebras at $q=0$}, 
J. Alg. Comb. {\bf 6} (1997), 339--376.
%
\bibitem{La} {\sc A.~Lascoux}, 
{\it Anneau de Grothendieck de la vari\'et\'e de drapeaux}, 
in The Grothendieck Festschrift, Vol. III,,
Progr. Math., {\bf 88} (1990), 1--34.
%
%
\bibitem{Loth}  {\sc A.~Lascoux, B.~Leclerc} and {\sc J.-Y.~Thibon}, 
{\it The plactic monoid}, Chapter 5 of {\sc M. Lothaire}, 
{\it Algebraic Combinatorics on Words}, 
Cambridge University Press, to appear.
%
\bibitem{LR}{\sc J.-L. Loday} and {\sc M. Ronco}, 
{\it Hopf algebra of the planar binary trees}, 
Adv. in Math. {\bf 139} (1998), 293--309.      
%
%
%
\bibitem{Mv}{\sc C. Malvenuto}, {\it $P$-partitions and the plactic congruence},
 Graphs Combin. {\bf 9} (1993),  63--73.
%
\bibitem{MR1} {\sc C.~Malvenuto} and {\sc C.~Reutenauer}, 
{\it Duality between quasi-symmetric functions and Solomon descent
algebra}, 
J. Algebra, {\bf 177} (1995), 967--982.

\bibitem{MV} {\sc R. Mart\'\i nez-Villa}, 
{\it Applications of Koszul algebras: The preprojective algebra},
in Representation Theory of Algebras,
CMS Conference Procceedings, vol. 18 (1996), 487--504.
%
\bibitem{No}{\sc P.N. Norton}, 
{\it $0$-Hecke algebras},
J. Austral. Math. Soc. Ser. A {\bf 227} (1979), 337--357.
%
\bibitem{Nov}{\sc J.-C.~Novelli}, 
{\it On the hypoplactic monoid},
Proceedings of the 8th conference ``Formal Power Series and
Algebraic Combinatorics'', Vienna, 1997.
%
\bibitem{PR}{\sc S. Poirier} and {\sc C. Reutenauer}, 
{\it Alg\`ebre de Hopf des tableaux}, 
Ann. Sci. Math. Q\'ebec {\bf 19} (1995), 79--90.
%
\bibitem{Re}
{\sc C.~Reutenauer},
{\it Free Lie Algebras},
{Oxford science publications}, (1993).
%
\bibitem{St1}{\sc R.P. Stanley}, 
{\it Ordered structures and partitions}, 
Memoirs Amer. Math. Soc. {\bf 119} (1972).
%
\bibitem{St}{\sc R.P. Stanley}, 
{\it Generalized riffle shuffles and quasisymmetric functions}, 
preprint, 1999.
%
\bibitem{TU}{\sc J.-Y.~Thibon} and {\sc B.-C.-V.~Ung}, 
{\it Quantum quasi-symmetric functions and Hecke algebras}, 
J. Phys. A: Math. Gen., {\bf 29} (1996),
7337--7348.
%
\bibitem{U}{\sc B.-C.-V. Ung},
{\it Combinatorial identities for series of quasi-symmetric
functions}, 
preprint IGM 98-15, Universit\'e de Marne-la-Vall\'ee.
\end{thebibliography}
\end{document}